\documentstyle[preprint,aps,prb]{revtex}
\def\C{{C\kern-.647em C}}
\def\MC{M_{\C}}

\def\MVC{M_{V,\C}}
\def\MVrC{M_{V^{r},\C}}
\def\MVpC{M_{V^{p},\C}}

\def\sqtimes{\Box{\raise 0.1em\hbox{\kern-0.76em $\times$}}}
\def\twsqtimes{\widetilde{\raise-0.1em\hbox{$\sqtimes$}}}

\begin{document}

\title{Abstract carrier space formalism for\\ the irreducible tensor 
operators of\\ compact quantum
group algebras}
\author{J.F.Cornwell} 
\address{School of Physics and Astronomy, University of 
St.Andrews,\\North Haugh, St.Andrews, Fife,
KY16 9SS, Scotland, U.K.} 
\date{October 12, 1995} 
\maketitle
\begin{abstract} Defining conditions for irreducible tensor operators 
associated with the
unitary irreducible corepresentations of compact quantum group algebras 
are deduced within the
framework of the abstract carrier space formalism. It is shown that there 
are {\em{two}} types of
irreducible tensor operator, which may be called `ordinary' and `twisted'. 
The consistency of the
definitions is demonstrated, and various consequences are deduced, 
including generalizations of the
Wigner-Eckart theorem for both the ordinary and twisted operators. 
Examples of irreducible tensor
operators for the standard deformation of the function algebra of the 
compact Lie group $SU(2)$ are
described to demonstrate the applicability of the new definitions.
\end{abstract}
\vskip 0.15in \centerline{{\sl Running title:} Irreducible tensor operators}
\newpage
\section{Introduction}

Most of the applications to physics of the theories of groups and Lie 
algebras
depend on the Wigner-Eckart theorem, and so it is of great interest to see 
how this theorem
generalizes to other algebraic structures. In a previous paper$^{1}$, 
hereafter referred to as Paper
I, a study of the definitions and properties of  irreducible tensor operators 
for a {\em{compact
quantum group algebra}} $\cal{A}$ was initiated by examining the case of 
the right regular and left
regular coaction formalisms, and was extended to the case of operators 
associated with
the corresponding quantum homogeneous spaces of $\cal{A}$. In the 
present paper this will be further extended to the case of
operators acting in the abstract carrier spaces of irreducible 
corepresentations of $\cal{A}$.

 The plan of the present paper is as follows. The remainder of this section 
will be devoted to
putting the analysis that follows into context, first by reviewing briefly the 
situation for compact
Lie groups, and then by indicating the background for the generalization 
of compact Lie groups to
compact quantum group algebras. In the next section the most relevant 
properties of compact quantum
group algebras will be summarized, particular attention being devoted to 
the essential role played by
corepresentations. This summary is continued in Section~III, where the 
two
different types of tensor product of corepresentations and their associated 
Clebsch-Gordan
coefficients are briefly discussed.  The heart of the paper is reached in  
Section~IV, where the
irreducible tensor operators are defined and some of their immediate 
properties are deduced. In
particular, it will be shown there that  there are {\em{two}} types of 
irreducible tensor operators,
which will be described as being {\em{ordinary}} and {\em{twisted}} 
respectively. The
{\em{motivations}} for the definitions of Section~IV are deliberately
relegated to  the Appendix in order to emphasize that the treatment given 
for the compact quantum
group algebras in Sections~IV and
V entirely self-contained.  In Section~V it is shown
that there are {\em{two}} theorems of the Wigner-Eckart type, one for the 
`ordinary' and one for the
`twisted' irreducible tensor operators. To illustrate the applicability of the 
new definitions of
Section~IV, examples of irreducible tensor operators for the standard
deformation of the function algebra of the compact Lie group $SU(2)$ are 
described in
Section~VI. Unless otherwise stated, the notations, definitions, and 
terminology are exactly the same
as those given in Paper I, which also contains an account of the 
relationship of the present line of
study to previous work on the applicability of the Wigner-Eckart theorem 
to quantum groups.

Because the space of functions defined on a compact Lie group $\cal{G}$ is 
a special example of a
compact quantum group algebra, all the well-known results for compact 
Lie groups naturally reappear
in this particular case. However, as the detailed analysis shows, the 
theory in the general
situation is rather more subtle, and exhibits various complications. 
Nevertheless, the point of view
of the present communication is best introduced by considering first the 
abstract carrier space
formalism in this very well established and  familiar context of a compact 
Lie group
$\cal{G}$ (c.f. Refs. 2, 3).  Let $V^{p}$ be a carrier space for a unitary 
irreducible
representation ${\bf{\Gamma}}^{p}$ of $\cal{G}$, and let $\psi^{p}_{1},
\psi^{p}_{2},\ldots, \psi^{p}_{d_{p}}$ be an ortho-normal basis for $V^{p}$. 
Define for each $T
\in \cal{G}$ a linear operator $\Phi^{p}(T)$ that acts on $V^{p}$ by the 
requirement that
\begin{equation} \Phi^{p}(T)(\psi^{p}_{n}) = \sum_{m=1}^{d_{p}}
\Gamma^{p}(T)_{mn}\psi^{p}_{m}\label{eq:Vintro1}\end{equation}  
for all $T \in \cal{G}$ and all $n =
1, 2, \ldots, d_{p}$. Let ${\bf{\Gamma}}^{p}$,
${\bf{\Gamma}}^{q}$, and ${\bf{\Gamma}}^{r}$ be any three unitary 
irreducible resentations of
$\cal{G}$. Then one can consider a set of irreducible tensor operators 
$Q^{q}_{1}, Q^{q}_{2},
\ldots, Q^{q}_{d_{q}}$ that each map
$V^{p}$ into $V^{r}$ and which are such that
\begin{equation}\Phi^{r}(T)\,Q^{q}_{n}\,\Phi^{p}(T)^{-1} = 
\sum_{m=1}^{d_{q}}
\Gamma^{q}(T)_{mn}Q^{q}_{m}\label{eq:Vintro2}\end{equation} 
for all $T \in \cal{G}$ and all $n = 1,
2, \ldots, d_{q}$. In this case the Wigner-Eckart theorem deals with inner 
products $\langle\;,\;
\rangle$ defined on
$V^{r}$ and states that the $j, k,$ and $\ell$ dependence of 
$\langle \psi^{r}_{\ell},Q^{q}_{k}(\psi^{p}_{j}) \rangle$  depends only on  
Clebsch-Gordan
coefficients for the reduction of the tensor product ${\bf{\Gamma}}^{p} 
\otimes {\bf{\Gamma}}^{q}$
into its irreducible constituents ${\bf{\Gamma}}^{r}$. 

In a minor extension of this formalism, one
could introduce an inner product space $V$ that is a direct sum of carrier 
spaces of certain unitary
irreducible representations of $\cal{G}$ and which contains at least 
$V^{p}\oplus V^{r}$ (and which,
in the extreme case,  may contain one carrier space for every inequivalent 
irreducible
representation of $\cal{G}$).  Then, for each $T \in \cal{G}$ an operator 
$\Phi(T)$ can be defined
which maps elements of $V$ into $V$, and which acts as
$\Phi^{p}(T)$ on $V^{p}$, as $\Phi^{r}(T)$ on $V^{r}$, and so on.  The 
irreducible tensor operators
are then required to each map
$V$ into $V$ and to be such that
$\Phi(T)\,Q^{q}_{n}\;\Phi(T)^{-1} = \sum_{m=1}^{d_{q}}
\Gamma^{q}(T)_{mn}Q^{q}_{m}$ for all $T \in \cal{G}$ and all $n = 1, 2, 
\ldots, d_{q}$. In this case
the Wigner-Eckart theorem deals with inner products $\langle\;,\; 
\rangle$ defined on
$V$, but is otherwise the same as above. 

As emphasized in Paper I, one most important lesson that can be drawn 
from these simple group
theoretical considerations concerns the {\em{consistency}} of the 
definitions of the basis vectors
{\em{and}} of the irreducible tensor operators.   As
$\Phi^{p}(T)\Phi^{p}(T^{\prime}) = \Phi^{p}(TT^{\prime})$ and
${\bf{\Gamma}}^{p}(T){\bf{\Gamma}}^{p}(T^{\prime}) = 
{\bf{\Gamma}}^{p}(TT^{\prime})$ for all
$T,T^{\prime} \in \cal{G}$, it follows that if (\ref{eq:Vintro1}) is valid for 
$T$ and for
$T^{\prime}$, then it is also valid for their product
$TT^{\prime}$. Similarly, and very significantly, by defining for each $T 
\in \cal{G}$ an operator
$\Psi(T)$ by 
\begin{equation} \Psi(T)(Q) = \Phi^{r}(T)\,Q\,\Phi^{p}(T)^{-
1}\label{eq:Vext1}\end{equation} 
for every
operator
$Q$ that maps
$V^{p}$ into $V^{r}$, the definition (\ref{eq:Vintro2}) can be recast as  
\begin{equation} \Psi(T)(Q^{q}_{n}) = \sum_{m=1}^{d_{q}}
\Gamma^{q}(T)_{mn}Q^{q}_{m} \label{eq:5.34a} \end{equation}      for all 
$T \in \cal{G}$ and all $n
= 1, 2, \ldots, d_{q}$. As $\Psi(T)\Psi(T^{\prime}) = \Psi(TT^{\prime})$  
for all
$T,T^{\prime} \in \cal{G}$, it follows that if (\ref{eq:5.34a}) is valid for 
$T$ and for
$T^{\prime}$, then it is also valid for their product $TT^{\prime}$. Put 
another way, because of the
similarity in form between (\ref{eq:Vintro1}) and (\ref{eq:5.34a}), the 
{\em{consistency}} of the
definition (\ref{eq:Vintro2}) of the irreducible tensor operators $Q^{q}_{n}$ 
is ensured by the fact
that they too form a basis for a carrier space, this time for  
${\bf{\Gamma}}^{q}$. In the
analysis that follows (cf. Section~IV), essentially this argument will be
used to justify the definitions that will be given for the irreducible tensor 
operators of the
compact quantum group algebras, the only essential
difference being that the argument has to be cast in terms of 
{\em{co}}representations instead of
representations.

As is well known, the set of functions defined on a Lie group $\cal{G}$ 
form a Hopf algebra,
$\cal{A}$, and  the dual ${\cal{A}}^{\prime}$ of $\cal{A}$ is the universal 
enveloping algebra
of the Lie algebra $\cal{L}$ of $\cal{G}$. Moreover, the structure of 
$\cal{G}$ can be encoded into
the structure of $\cal{A}$, and, in particular,
$\cal{A}$ is commutative.  A `deformation' (or `quantization') of  
${\cal{A}}^{\prime}$ induces a
corresponding deformation of $\cal{A}$, and will make $\cal{A}$ non-
commutative as well as being
non-cocommutative. Although most attention has been
focused on the deformed Hopf algebras ${\cal{A}}^{\prime}$, it has been 
demonstrated by the
pioneering work of Woronowicz$^{4-6}$, which itself has been refined and 
developed by
Dijkhuizen and Koornwinder$^{7-12}$, that it is of very great interest to 
produce a
self-contained and direct study of generalizations of the Hopf algebras
$\cal{A}$. This can be done by assuming that they have certain 
characteristic properties, and the
resulting structures have been  called {\em{compact matrix 
pseudogroups}} by
Woronowicz$^{4-6}$, and {\em{compact quantum group algebras}} by 
Dijkhuizen and
Koornwinder$^{7-12}$. It is these that provide the framework for the  
present paper, which,
as intimated above, is devoted to the study of the irreducible tensor 
operators for compact quantum
group algebras in the abstract carrier space formalism.

\section{Corepresentations of compact quantum group algebras}

It should be recalled (c.f. Refs. 1, 7-12) that a {\em{right $\cal{A}$-
comodule}} consists of a
vector space $V$ and a linear mapping $\pi_{V}$ from
$V$ to $V \otimes \cal{A}$ such that
\begin{equation}  (\pi_{V} \otimes id) \circ \pi_{V} = (id \otimes \Delta) 
\circ \pi_{V}
\label{eq:R1} \end{equation} and \begin{equation}   ((id \otimes \epsilon) 
\circ \pi_{V})(v) =
v \otimes 1_{\C}
\label{eq:R2} \end{equation} 
for all $v \in V$. The operation $\pi_{V}$ is then said to be a {\em{right 
coaction}} and provides a
{\em{corepresentation}} of $\cal{A}$ with carrier space $V$.  

Of great importance are the
finite-dimensional {\em{irreducible}} corepresentations, which, for a 
compact quantum group algebra
$\cal{A}$, are assumed to form a {\em{countable}} set (up to equivalence). 
Moreover each
such irreducible corepresentation is equivalent to a {\em{unitary}} 
corepresentation. Let $\pi^{p}$,
for
$p = 1, 2,
\ldots$, denote the set of unitary  irreducible corepresentations of 
$\cal{A}$ (one being chosen from
every equivalence class), and let 
$V^{p}$ be a carrier space of $\pi^{p}$, assumed to be of finite dimension
$d_{p}$, with basis $v^{p}_{1}, v^{p}_{2}, \ldots, v^{p}_{d_{p}}$. Then there 
exists a uniquely
determined set of elements
$\pi^{p}_{jk}$ of $\cal{A}$ (for
$j,k = 1, 2, \ldots, d_{p}$), called the {\em{matrix coefficients}} of 
$\pi^{p}$, which are such that
\begin{equation}  \pi^{p}(v^{p}_{j}) = \sum_{k=1}^{d_{p}} v^{p}_{k} 
\otimes \pi^{p}_{kj}
\label{eq:R5} \end{equation} for all $j = 1, 2, \ldots, d_{p}$.  The 
requirements (\ref{eq:R1}) and
(\ref{eq:R2}) then imply that
\begin{equation}  \Delta(\pi^{p}_{jk}) = \sum_{\ell=1}^{d_{p}} 
\pi^{p}_{j\ell} \otimes \pi^{p}_{\ell
k}
\label{eq:R6} \end{equation} and
\begin{equation}  \epsilon(\pi^{p}_{jk}) = \delta_{jk}
\label{eq:R7} \end{equation} (for $j,k = 1, 2,\ldots, d_{p}$). The unitary 
requirement on $\pi^{p}$
implies that
\begin{equation} S(\pi^{p}_{jk}) = \pi^{p*}_{kj} ,
\label{eq:D.sect1.1.45ii} \end{equation}
\begin{equation} \sum_{\ell=1}^{d_{p}} M(\pi^{p*}_{\ell j} \otimes 
\pi^{p}_{\ell k}) = \delta
_{jk} 1_{\cal{A}}  ,
\label{eq:D.sect1.1.45iii} \end{equation} and
\begin{equation} \sum_{\ell=1}^{d_{p}} M(\pi^{p}_{j\ell} \otimes 
\pi^{p*}_{k\ell}) = \delta
_{jk} 1_{\cal{A}}  
\label{eq:D.sect1.1.45iv} \end{equation} (for all $j,k = 1, 2, \ldots, d_{p}$). 
For a compact
quantum group algebra $\cal{A}$ this set of matrix coefficients are 
assumed (c.f.Refs. 1, 7-12) to
form a basis for
$\cal{A}$.

Let ${\cal{P}}^{p}_{ij}$ be a set of projection operators for $V^{p}$ that are 
defined
by
\begin{equation} {\cal{P}}^{p}_{ij}(v^{p}_{k}) = \delta_{ik}\; 
v^{p}_{j}\label{eq:V3}\end{equation}
for all $i,j,k = 1, 2, \ldots, d_{p}$. Let $\pi^{r}$ be another unitary  
irreducible corepresentation
of $\cal{A}$, and let ${\cal{L}}^{pr}$ be the set of linear operators that 
map elements of $V^{p}$
into
$V^{r}$. A basis for ${\cal{L}}^{pr}$ is provided by the set of operators 
${\cal{P}}^{pr}_{ij}$ that
are defined by
\begin{equation} {\cal{P}}^{pr}_{ij}(v^{p}_{k}) = \delta_{ik}\; 
v^{r}_{j}\label{eq:V5}\end{equation}
for all $i,k = 1, 2, \ldots, d_{p}$ and all $j = 1, 2, \ldots, d_{r}$. Then
\begin{equation} {\cal{P}}^{r}_{mn} \circ {\cal{P}}^{pr}_{k\ell}\circ 
{\cal{P}}^{p}_{ij} =
\delta_{kj}\; \delta_{m\ell}\; 
{\cal{P}}^{pr}_{in}\label{eq:V9}\end{equation} for all $i,j,k =
1, 2, \ldots, d_{p}$ and all $\ell,m,n = 1, 2, \ldots, d_{r}$. If $Q$ is any 
element of
${\cal{L}}^{pr}$, then
\begin{equation} Q(v^{p}_{k}) = \sum_{j = 1}^{d_{r}} q_{jk}\; 
v^{r}_{j}\label{eq:V8}\end{equation}
for all $k = 1, 2, \ldots, d_{p}$, where $q_{jk}$ are the complex numbers 
that are defined by
\begin{equation} q_{jk} = \langle v^{r}_{j}, 
Q(v^{p}_{k})\rangle\label{eq:V7}\end{equation}
for all $k = 1, 2, \ldots, d_{p}$ and all $j = 1, 2, \ldots, d_{r}$, 
$\langle\;,\;\rangle$ being the
inner product of $V^{r}$. Moreover one can write
\begin{equation} Q = \sum_{i=1}^{d_{p}} \sum_{j=1}^{d_{r}} 
\;q_{ji}\,{\cal{P}}^{pr}_{ij}\;.
\label{eq:V6}\end{equation}

\section{Tensor products and Clebsch-Gordan coefficients}

\subsection{Ordinary and twisted tensor products}

With the {\em{ordinary tensor product}} of two irreducible 
corepresentations
$\pi^{p}$ and $\pi^{q}$ of $\cal{A}$ (with carrier spaces $V^{p}$ and 
$V^{q}$ respectively) being
defined as the mapping
$\pi^{p} \sqtimes \pi^{q}$ from $V^{p} \otimes V^{q}$ to $V^{p} \otimes 
V^{q} \otimes \cal{A}$ that
is such that
\begin{equation} \pi^{p} \sqtimes \pi^{q} = (id \otimes id \otimes M) 
\circ
(id \otimes \sigma \otimes id) \circ (\pi^{p} \otimes \pi^{q}) \; ,
\label{eq:V39} \end{equation}
it is easily shown from (\ref{eq:R5}) that the corresponding matrix 
coefficients are given by 
\begin{equation} (\pi^{p} \sqtimes \pi^{q})_{st,jk} = M(\pi_{sj}^{p} 
\otimes \pi_{tk}^{q}) 
\label{eq:V42} \end{equation}
for all $j,s = 1, 2, \ldots, d_{p}$ and all $k,t = 1, 2, \ldots, d_{q}$.

Similarly, with the {\em{twisted tensor product}} of $\pi^{p}$ and 
$\pi^{q}$  being defined as the
mapping
$\pi^{p} \twsqtimes \pi^{q}$ from $V^{p} \otimes V^{q}$ to $V^{p} 
\otimes V^{q} \otimes \cal{A}$ that
is such that
\begin{equation} \pi^{p} \twsqtimes \pi^{q} = (id \otimes id \otimes M) 
\circ
(id \otimes id \otimes \sigma) \circ (id \otimes \sigma \otimes id) \circ 
(\pi^{p} \otimes \pi^{q})
\; ,
\label{eq:V39a} \end{equation}
it is also easily shown from (\ref{eq:R5}) that the corresponding matrix 
coefficients are given by 
\begin{equation} (\pi^{p} \twsqtimes \pi^{q})_{st,jk} = M(\pi_{tk}^{q} 
\otimes \pi_{sj}^{p}) 
\label{eq:V42a} \end{equation}
for all $j,s = 1, 2, \ldots, d_{p}$ and all $k,t = 1, 2, \ldots, d_{q}$.

\subsection{Clebsch-Gordan coefficients}

Suppose that the ordinary tensor product $\pi^{p} \sqtimes \pi^{q}$ is 
reducible (and hence is
completely reducible (c.f. Refs. 1, 7-12)), and that
$n_{pq}^{r}$ is the number of times that the irreducible corepresentation
$\pi^{r}$ (or a corepresentation equivalent to it) appears in its reduction. 
If the carrier spaces
$V^{p}$ and $V^{q}$ have basis elements
$v^{p}_{1},v^{p}_{1}, \ldots, v^{p}_{d_{p}}$ and 
$v^{q}_{1},v^{q}_{1}, \ldots, v^{q}_{d_{q}}$ respectively, then the set of 
elements $v^{p}_{j}
\otimes v^{q}_{k}$ form a basis for $V^{p} \otimes V^{q}$, the carrier 
space of $\pi^{p} \sqtimes
\pi^{q}$, and consequently appropriate linear combinations form bases for 
all the irreducible
corepresentations
$\pi^{r}$ that appear in the reduction of the tensor product. Let
$w^{r,\alpha}_{\ell}$ be such a combination, so that 
\begin{equation} w^{r,\alpha}_{\ell} = \sum_{j=1}^{d_{p}} 
\sum_{k=1}^{d_{q}} \left( \begin{array}{cc}
p & q\\ j &  k \end{array} \right| \left. \begin{array}{ccc}r&,& 
\alpha\\\ell & &  
\end{array}\right) v^{p}_{j} \otimes v^{q}_{k} ,
\label{eq:C43} \end{equation} for $\ell = 1, 2, \ldots, d_{r}$, and $\alpha 
= 1, 2, \ldots,
n_{pq}^{r}$, and 
\begin{equation} (\pi^{p} \sqtimes \pi^{q})(w^{r,\alpha}_{\ell})  = 
\sum_{u=1}^{d_{r}}
w^{r,\alpha}_{u} \otimes \pi^{r}_{ul} ,
\label{eq:C44} \end{equation} for $u = 1, 2, \ldots, d_{r}$, and $\alpha = 
1, 2, \ldots,
n_{pq}^{r}$. The inverse of (\ref{eq:C43}) is
\begin{equation}  v^{p}_{j} \otimes v^{q}_{k} = \sum_{r} 
\sum_{\alpha=1}^{n_{pq}^{r}} \sum_{\ell
=1}^{d_{r}} \left(
\begin{array}{ccc} r&,& \alpha\\\ell & &     \end{array} \right| \left. 
\begin{array}{cc}p & q\\ j &
k   \end{array}\right) w^{r,\alpha}_{\ell} ,
\label{eq:C45} \end{equation} for $j = 1, 2, \ldots, d_{p}$ and $k = 1, 2, 
\ldots, d_{q}$.  The
{\em{Clebsch-Gordan coefficients}} defined in (\ref{eq:C43}) form the 
elements of a $d_{p} \times
d_{q}$ matrix ${\bf{C}}$, while the inverse coefficients defined in 
(\ref{eq:C45}) form the elements
of 
${\bf{C}}^{-1}$, where
\begin{equation}{\bf{C}}^{-1} ({\bf{\pi}}^{p} \sqtimes 
{\bf{\pi}}^{q}){\bf{C}} = \sum_{r} \oplus
n_{pq}^{r} {\bf{\pi}}^{r} .
\label{eq:R159} \end{equation} This implies that
\begin{equation}  (\pi^{p} \sqtimes \pi^{q})_{is,jt} =  \sum_{r} 
\sum_{\alpha = 1}^{n_{pq}^{r}} \sum_{\ell,u=1}^{d_{r}} \left(
\begin{array}{cc} p & q\\ i &  s \end{array} \right| \left. 
\begin{array}{ccc}r&,& \alpha\\u & &  
\end{array}\right) \pi^{r}_{u\ell} \left(
\begin{array}{ccc}r&,& \alpha\\\ell & & \end{array} \right| \left.
\begin{array}{cc} p & q\\ j &  t  
\end{array}\right)
\label{eq:C45q} \end{equation} for $i,j = 1, 2, \ldots, d_{p}$, and $s,t = 1, 
2, \ldots, d_{q}$.

Thus, by (\ref{eq:V42}), the product of any two basis elements of $\cal{A}$ 
can be expressed in
terms of Clebsch-Gordan coefficients, for
\begin{equation}  M(\pi^{p}_{ij} \otimes \pi^{q}_{st}) =  \sum_{r} 
\sum_{\alpha =
1}^{n_{pq}^{r}}
\sum_{\ell,u=1}^{d_{r}} \left(
\begin{array}{cc} p & q\\ i &  s \end{array} \right| \left. 
\begin{array}{ccc}r&,& \alpha\\u & &  
\end{array}\right) \pi^{r}_{u\ell} \left(
\begin{array}{ccc}r&,& \alpha\\\ell & & \end{array} \right| \left.
\begin{array}{cc} p & q\\ j &  t  
\end{array}\right)
\label{eq:C45qextra} \end{equation} for $i,j = 1, 2, \ldots, d_{p}$, and $s,t 
= 1, 2, \ldots, d_{q}$.
This is essentially the converse of the relation (I.137), that is, of 
\begin{eqnarray} \begin{array}{lcl}h(\pi^{r*}_{ul}\pi^{q}_{tk}\pi^{p}_{sj})
 &= & \sum_{\alpha=1}^{n_{qp}^{r}} \sum_{v=1}^{d_{r}} 
\left(\begin{array}{ccc} r&,&
\alpha\\\ell & &    
\end{array} \right| \left. \begin{array}{cc}q & p\\ k & j   
\end{array}\right) \left(
\begin{array}{cc} q & p\\ t &  s \end{array} \right| \left. 
\begin{array}{ccc}r&,& \alpha\\v & &  
\end{array}\right)   \\
  & &\times  \{(({\bf{F}}^{r})^{-1})_{vu}/tr(({\bf{F}}^{r})^{-1})\}\end{array}
\label{eq:R164} \end{eqnarray} for all $j = 1, 2, \ldots, d_{p}$, $k = 1, 2, 
\ldots, d_{q}$, and $l
= 1, 2, \ldots, d_{r}$. Here $h$ is the Haar functional and ${\bf{F}}^{r}$ is 
a
non-singular  $d_{r} \times d_{r}$ matrix with the property that
\begin{equation} \sum_{k=1}^{d_{r}} F_{jk}^{r} \pi_{k\ell}^{r} = 
\sum_{k=1}^{d_{r}}
\pi_{jk}^{r\ddagger} F_{k\ell}^{r}
\label{eq:R148} \end{equation} (for all $j,\ell = 1,2,\ldots,d_{r}$), where 
$\pi^{r\ddagger}$
is the doubly contragredient partner of $\pi^{r}$.

\section{The irreducible tensor operators}

\subsection{Introduction}

Let $\pi^{p}$,$\pi^{q}$, and $\pi^{r}$ be unitary irreducible right 
coactions of $\cal{A}$ of
dimensions $d_{p}$, $d_{q}$, and $d_{r}$ respectively, and with matrix 
coefficients $\pi^{p}_{jk}$,
$\pi^{q}_{jk}$, and $\pi^{r}_{jk}$ respectively. Let ${\cal{L}}^{pr}$ be the 
vector space of
operators introduced in Section~II. It will be shown that  there exist two
types of irreducible tensor operators that are members of ${\cal{L}}^{pr}$ 
and which belong to the
corepresentation
$\pi^{q}$. These will be denoted by $Q^{q}_{j}$ and 
${\widetilde{Q}}^{q}_{j}$ (for $j = 1, 2, \ldots,
d_{q}$), and will be called {\em{ordinary}} and {\em{twisted}} irreducible 
tensor operators
respectively. Naturally the two types of irreducible tensor operators 
coincide in the special case
in which $\cal{A}$ is commutative.    

It will also be shown that the definitions of both of these types of 
irreducible tensor operators is
easily extended to the case in which $V$ is a vector space that is a direct 
sum of carrier
spaces of unitary irreducible corepresentations of $\cal{A}$ and which 
contains at least
$V^{p}\oplus V^{r}$. 

\subsection{Definitions of irreducible tensor operators}

\subsubsection{Definition of the ordinary irreducible tensor operators 
$Q^{q}_{j}$}

The {\em{ordinary irreducible tensor operators}}
$Q^{q}_{j}$ belonging to the unitary irreducible right coaction $\pi^{q}$ of
$\cal{A}$ are {\em{defined}} to be members of ${\cal{L}}^{pr}$ that satisfy 
the condition 
\begin{equation}  ((id \otimes M) \circ (\pi^{r} \otimes id) \circ (Q^{q}_{j}
\otimes S) \circ \pi^{p})(v^{p}) = \sum_{k=1}^{d_{q}} Q^{q}_{k}(v^{p}) 
\otimes
\pi^{q}_{kj}
\label{eq:V20} \end{equation}  for all $v^{p} \in V^{p}$ and all $j = 1, 2, 
\ldots, d_{q}$.
Clearly this definition involves {\em{only}} quantities defined for 
$\cal{A}$ and its coactions.
Both sides (\ref{eq:V20}) are members of
$V^{r} \otimes \cal{A}$. (The
motivation behind the definition (\ref{eq:V20}) is explained in  Section~B 
of the Appendix.)

It will now be shown that (\ref{eq:V20}) provides a {\em{consistent}} 
definition, in that it can be
re-expressed by saying that the operators $Q^{q}_{j}$ (for $j = 1, 2, \ldots, 
d_{q}$) form the
basis of an irreducible subspace of a carrier space for a certain right 
coaction of $\cal{A}$. This
right coaction will be denoted by $\pi_{{\cal{L}}^{pr}}$, as its carrier space 
is 
${\cal{L}}^{pr}$. The {\em{definition}} of $\pi_{{\cal{L}}^{pr}}$ is then 
that it is the mapping of 
$\pi_{{\cal{L}}^{pr}}$ into
$\pi_{{\cal{L}}^{pr}} \otimes \cal{A}$ that specified by
\begin{equation} \pi_{{\cal{L}}^{pr}}(Q) =
\sum_{i,j=1}^{d_{p}}\;\sum_{m,n=1}^{d_{r}}\; (q_{ni}\,{\cal{P}}^{pr}_{jm})
\otimes M(\pi^{r}_{mn} \otimes S(\pi^{p}_{ij}))\; ,
\label{eq:V30} \end{equation}  for all $Q \in {\cal{L}}^{pr}$, where 
$q_{ni}$ and
${\cal{P}}^{pr}_{jm}$ are defined in (\ref{eq:V7}) and (\ref{eq:V5}). (The
motivation for the definition (\ref{eq:V30}) is given in Section~B of the 
Appendix.)     

It is then quite easily shown that
$\pi_{{\cal{L}}^{pr}}$ satisfies (\ref{eq:R1}) and (\ref{eq:R2}) (with 
$\pi_{V}$ and $V$
replaced by $\pi_{{\cal{L}}^{pr}}$ and  ${\cal{L}}^{pr}$ respectively), and 
hence
$\pi_{{\cal{L}}^{pr}}$ is indeed a right coaction with carrier space
${\cal{L}}^{pr}$. Moreover, it is easily demonstrated that
\begin{equation} (\pi_{{\cal{L}}^{pr}}(Q))(v^{p} \otimes 1_{{\cal{A}}}) = 
((id \otimes M)  \circ
(\pi^{r} \otimes id) \circ (Q \otimes S) \circ \pi^{p})(v^{p})
\label{eq:V33} \end{equation}  
for all $v^{p} \in V^{p}$ and all $Q \in {\cal{L}}^{pr}$.
Thus (\ref{eq:V20}) and (\ref{eq:V33}) imply that the definition 
(\ref{eq:V20}) can be written
equivalently as
\begin{equation} \pi_{{\cal{L}}^{pr}}(Q^{q}_{j}) = \sum_{k=1}^{d_{q}} 
Q^{q}_{k} \otimes
\pi^{q}_{kj}  
\label{eq:V34} \end{equation}  (for all $j = 1, 2, \ldots, d_{q}$). Because 
(\ref{eq:V34}) is
similar in form to (\ref{eq:R5}), and as
$\pi_{{\cal{L}}^{pr}}$ is a right coaction with carrier space 
${\cal{L}}^{pr}$, the
consistency of the definition (\ref{eq:V20}) is now ensured.

Now consider the situation in which $V$ is a vector space that is a direct 
sum of carrier
spaces of unitary irreducible corepresentations of $\cal{A}$ and which 
contains at least
$V^{p}\oplus V^{r}$. Let $\pi$ be the mapping of $V$ into $V \otimes 
\cal{A}$ that coincides
with $\pi^{p}$ on $V^{p}$ and with $\pi^{r}$ on $V^{r}$, and which acts 
similarly on any other
carrier spaces that might be contained in $V$. Then the generalization of 
(\ref{eq:V20}) is clearly 
\begin{equation}  ((id \otimes M) \circ (\pi \otimes id) \circ (Q^{q}_{j}
\otimes S) \circ \pi)(v) = \sum_{k=1}^{d_{q}} Q^{q}_{k}(v) \otimes
\pi^{q}_{kj}
\label{eq:V20ext} \end{equation}  for all $v \in V$ and all $j = 1, 2, \ldots, 
d_{q}$. (The
consistency of the definition (\ref{eq:V20ext}) is an immediate 
consequence of the consistency of
(\ref{eq:V20})).

\subsubsection{Definition of the twisted irreducible tensor operators
${\widetilde{Q}}^{q}_{j}$}

The {\em{twisted irreducible tensor operators}} ${\widetilde{Q}}^{q}_{j}$ 
belonging to the unitary
irreducible right coaction $\pi^{q}$ of $\cal{A}$ are {\em{defined}} to be 
members of
${\cal{L}}^{pr}$ that satisfy the condition 
\begin{equation}  ((id \otimes M) \circ (id \otimes \sigma) \circ (\pi^{r} 
\otimes id) \circ
({\widetilde{Q}}^{qR}_{j}
\otimes S^{-1}) \circ \pi^{p})(v^{p}) = \sum_{k=1}^{d_{q}} 
{\widetilde{Q}}^{q}_{k}(v^{p}) \otimes
\pi^{q}_{kj}
\label{eq:V37e} \end{equation}  for all $v^{p} \in V^{p}$ and all $j = 1, 2, 
\ldots, d_{q}$. This
definition (\ref{eq:V37e})  differs from the corresponding definition
(\ref{eq:V20})  only in the replacement of $M$ by $M \circ \sigma$ 
{\em{and}} $S$
by $S^{-1}$ (neither of which have any effect in the special case in which 
$\cal{A}$ is
commutative). (See Section~B of the Appendix for further
discussion of this pair of substitutions. It should be recorded that 
Rittenberg and Scheunert$^{13}$
noted previously, in the context of what was essentially the abstract 
carrier space formalism
 of Section~I as generalized to irreducible {\em{representations}}  of the
dual 
${\cal{A}}^{\prime}$,  that these substitutions do produce another type of 
irreducible tensor
operator, but they did not pursue  this observation.)  

The demonstration that (\ref{eq:V37e}) provides a {\em{consistent}} 
definition again involves
showing  that it can be re-expressed by saying that the operators 
${\widetilde{Q}}^{q}_{j}$ (for $j
= 1, 2,
\ldots, d_{q}$) form the basis of an irreducible subspace of a carrier space 
for another  right
coaction
$\widetilde{\pi}_{{\cal{L}}^{pr}}$ of
$\cal{A}$. This right coaction is {\em{defined}} as the mapping of  
${\cal{L}}^{pr}$ into
${\cal{L}}^{pr} \otimes \cal{A}$ that specified by
\begin{equation} \widetilde{\pi}_{{\cal{L}}^{pr}}(Q) = 
\sum_{i,j=1}^{d_{p}}\;\sum_{m,n=1}^{d_{r}}\; (q_{ni}\,{\cal{P}}^{pr}_{jm})
\otimes M(S^{-1}(\pi^{p}_{ij} \otimes \pi^{r}_{mn}))\; 
\label{eq:V37f} \end{equation} 
for all $Q \in {\cal{L}}^{pr}$. (The motivation for the
definition (\ref{eq:V37f}) is given in Section~C of the
Appendix).     Then
\begin{equation} (\widetilde{\pi}_{{\cal{L}}^{pr}}(Q))(v^{p} \otimes 
1_{{\cal{A}}}) =
((id \otimes M) \circ (id \otimes \sigma) \circ (\pi^{r} \otimes id) \circ
(Q \otimes S^{-1}) \circ \pi^{p})(v^{p})
\label{eq:V37g} \end{equation}  for all $v^{p} \in V^{p}$ and all $Q \in 
{\cal{L}}^{pr}$.
Thus (\ref{eq:V37e}) and (\ref{eq:V37g}) imply that the definition 
(\ref{eq:V37e}) can be written
equivalently as
\begin{equation} \widetilde{\pi}_{{\cal{L}}^{pr}}({\widetilde{Q}}^{qR}_{j}) 
=
\sum_{k=1}^{d_{q}} {\widetilde{Q}}^{qR}_{k} \otimes
\pi^{q}_{kj}  
\label{eq:V37h} \end{equation}  (for all $j = 1, 2, \ldots, d_{q}$), which 
then ensures its
consistency.

In the situation in which $V$ is a vector space that is a direct sum of 
carrier
spaces of unitary irreducible corepresentations of $\cal{A}$ and which 
contains at least
$V^{p}\oplus V^{r}$, and with the mapping $\pi$ from $V$ into $V 
\otimes \cal{A}$ that is  defined
in the end of the previous subsubsection, the generalization of 
(\ref{eq:V37e}) is clearly 
\begin{equation}  ((id \otimes M) \circ (id \otimes \sigma) \circ (\pi 
\otimes id) \circ
({\widetilde{Q}}^{qR}_{j}
\otimes S^{-1}) \circ \pi)(v) = \sum_{k=1}^{d_{q}} 
{\widetilde{Q}}^{q}_{k}(v) \otimes
\pi^{q}_{kj}
\label{eq:V37eext} \end{equation}  for all $v \in V$ and all $j = 1, 2, 
\ldots, d_{q}$. (Again, the
consistency of the definition (\ref{eq:V37eext}) is an immediate 
consequence of the consistency of
(\ref{eq:V37e})).

\subsection{Properties of irreducible tensor operators}

\subsubsection{The identity operator as an irreducible tensor operator}

Suppose that $V$ is a vector space that is a direct sum of carrier
spaces of unitary irreducible corepresentations of $\cal{A}$ and which 
contains at least
$V^{p}\oplus V^{r}$, and that $\pi$ is the mapping of $V$ into $V 
\otimes \cal{A}$ that is defined
in the previous subsection. Suppose that $Q$ is the {\em{identity 
operator}} $id$ of
$V$ (so that $Q(v) = v$ for all $v \in V$). Then, on using (\ref{eq:R1}) and 
(\ref{eq:R2}),
together with the Hopf algebra properties $M \circ (id \otimes S) \circ 
\Delta = u \circ \epsilon$
and $u(1_{\C}) = 1_{\cal{A}}$, it follows that
\begin{equation}  ((id \otimes M) \circ (\pi \otimes id) \circ (id
\otimes S) \circ \pi)(v) = v \otimes 1_{\cal{A}} 
\label{eq:V37cprime} \end{equation}  for all $v \in V$,  which, by 
(\ref{eq:V20ext}), leads to the
 conclusion that the identity operator $id$ is an {\em{ordinary}} 
irreducible tensor
operator for the one-dimensional {\em{identity}} corepresentation whose 
sole matrix coefficient is
$1_{{\cal{A}}}$.
 
It is easily checked (using (\ref{eq:V37eext})  in place of
(\ref{eq:V20ext})), that
$id$ is also a {\em{twisted}} irreducible tensor operator for this identity 
corepresentation.

The same conclusions for identity operators follow directly from 
(\ref{eq:V20}) and (\ref{eq:V37e})
in the special case in which $p = r$.

\subsubsection{Two useful identities for the ordinary irreducible tensor 
operators
$Q^{q}_{j}$ and 
${\widetilde{Q}}^{q}_{j}$}

 If  $Q^{q}_{k}$ is an {\em{ordinary}} irreducible tensor operator belonging 
to the unitary
irreducible corepresentation $\pi^{q}$ of $\cal{A}$ (as defined in 
(\ref{eq:V20})), and 
$v^{p}_{j}$ (for $j = 1, 2, \ldots, d_{p}$) provides a basis for the carrier 
space $V^{p}$ of the
unitary irreducible corepresentation $\pi^{p}$ of $\cal{A}$, then
\begin{equation} \pi^{r}(Q^{q}_{k}(v^{p}_{j})) =
\sum_{s=1}^{d_{p}} \sum_{t=1}^{d_{q}} (Q^{q}_{t}(v^{p}_{s}))
\otimes (M(\pi^{q}_{tk} \otimes \pi^{p}_{sj})) ,
\label{eq:V35} \end{equation}
for all $j = 1,2, \ldots, d_{p}$, and $k = 1,2, \ldots, d_{q}$.  By contrast, if 
${\widetilde{Q}}^{q}_{k}$ is a {\em{twisted}} irreducible tensor operator 
belonging
$\pi^{q}$, then
\begin{equation} \pi^{r}({\widetilde{Q}}^{q}_{k}(v^{p}_{j})) =
\sum_{s=1}^{d_{p}} \sum_{t=1}^{d_{q}} ({\widetilde{Q}}^{q}_{t}(v^{p}_{s}))
\otimes (M(\pi^{p}_{sj} \otimes \pi^{q}_{tk} )) ,
\label{eq:V37i} \end{equation}
for all $j = 1, 2, \ldots, d_{p}$ and $k = 1, 2, \ldots, d_{q}$. It should be 
noted that the factors
in the second term of the right-hand side of (\ref{eq:V37i}) are 
interchanged relative to those
of (\ref{eq:V35}).

The proof of (\ref{eq:V35}) is as follows. On applying (\ref{eq:R5}) and the
relation $S(\pi^{p}_{ks}) = \pi^{p*}_{sk}$, the left-hand side of  
(\ref{eq:V20}) (with $v^{p} =
v^{p}_{s}$) becomes
\[((id \otimes M) \circ (\pi^{r} \otimes id))(\sum_{k = 
1}^{d_{p}}(Q^{q}_{j}(v^{p}_{k})) \otimes
\pi^{p*}_{sk})\;\] 
On multiplying from the right with $id \otimes \pi^{p}_{si}$, summing 
over $s$,
and applying the relation $M(\pi^{p*}_{sk} \otimes \pi^{p}_{si}) = 
\delta_{ik}1_{\cal{A}}$, this
reduces to $\pi^{r}(Q^{q}_{i}(v^{p}_{j}))$. However,   multiplication of the 
right-hand side of
(\ref{eq:V20}) from the right with $id
\otimes \pi^{p}_{si}$  and summing over $s$ produces 
$\sum_{s=1}^{d_{p}} \sum_{t=1}^{d_{q}}
(Q^{q}_{t}(v^{p}_{s})) \otimes (M(\pi^{q}_{tk} \otimes \pi^{p}_{sj}))$. The 
line of proof for
(\ref{eq:V37i}) is similar.

\subsubsection{Identification of the corepresentations 
$\pi_{{\cal{L}}^{pr}}$  and
$\widetilde{\pi}_{{\cal{L}}^{pr}}$ of $\cal{A}$}

It is easily shown from the definition (\ref{eq:V30}) of 
$\pi_{{\cal{L}}^{pr}}$ that 
\begin{equation} \pi_{{\cal{L}}^{pr}}({\cal{P}}^{pr}_{ij}) =
\sum_{m=1}^{d_{p}}\sum_{n=1}^{d_{r}}{\cal{P}}^{pr}_{mn} \otimes 
(\pi^{r} \sqtimes
\bar{\pi}^{p})_{nm,ji}
\;,\label{eq:V45a}\end{equation} 
where ${\cal{P}}^{pr}_{ij}$ are the operators defined in (\ref{eq:V5}), and 
where $\bar{\pi}^{p}$ 
is the corepresentation of $\cal{A}$ that is {\em{conjugate}} to $\pi^{p}$, 
so that
its matrix coefficients are given by $\bar{\pi}_{jk}^{p} = 
(\pi^{p}_{jk})^{*}$. Then,
by (\ref{eq:V42}) and (\ref{eq:V42a}),
\begin{equation} \pi_{{\cal{L}}^{pr}}({\cal{P}}^{pr}_{ij}) =
\sum_{m=1}^{d_{p}}\sum_{n=1}^{d_{r}}{\cal{P}}^{pr}_{mn} \otimes 
(\bar{\pi}^{p} \twsqtimes
\pi^{r})_{mn,ij}
\;.\label{eq:V45c}\end{equation} 
This shows that $\pi_{{\cal{L}}^{pr}}$ is actually the right coaction that is 
given by the
{\em{twisted}} tensor product $\bar{\pi}^{p} \twsqtimes \pi^{r}$, and 
that the operators
${\cal{P}}^{pr}_{ij}$ are the basis vectors of the carrier space 
${\cal{L}}^{pr}$ of this coaction.

Taken with (\ref{eq:V34}), this indicates that the irreducible tensor 
operators $Q^{q}_{j}$ of
the definition (\ref{eq:V20}) exist only if $\pi^{q}$ is contained in the 
reduction of 
$\bar{\pi}^{p} \twsqtimes \pi^{r}$. As $\bar{\pi}^{p} \twsqtimes 
\pi^{r}$ and $\pi^{r}
\sqtimes \bar{\pi}^{p}$ are equivalent$^{1}$, this implies that 
$Q^{q}_{j}$ exists only if
$n^{q}_{r\bar{p}} > 0$. But $n^{q}_{r\bar{p}} = n^{r}_{qp}$ (c.f. (I.125)), so 
$Q^{q}_{j}$ exists
only if $n^{r}_{qp} > 0$. (These observations are confirmed by the explicit 
expressions for the
irreducible tensor operators given in (\ref{eq:V45d}) below and by the 
Wigner-Eckart theorem of
(\ref{eq:V53}) below).

Similarly, one can show from the definition (\ref{eq:V37f}) of 
$\widetilde{\pi}_{{\cal{L}}^{pr}}$
that 
\begin{equation} \widetilde{\pi}_{{\cal{L}}^{pr}}({\cal{P}}^{pr}_{ij}) =
\sum_{m=1}^{d_{p}}\sum_{n=1}^{d_{r}}{\cal{P}}^{pr}_{mn} \otimes 
(\overline{\pi^{p\ddagger}} \sqtimes
\pi^{r})_{mn,ij}
\;,\label{eq:V45c2}\end{equation} 
where $\overline{\pi^{p\ddagger}}$  is the corepresentation of $\cal{A}$ 
that is {\em{conjugate}} to
the corepresentation $\pi^{p\ddagger}$ that is itself {\em{doubly 
contragredient}} to $\pi^{p}$.
This shows that $\widetilde{\pi}_{{\cal{L}}^{pr}}$ is the right coaction 
that is given by
the ordinary tensor product $\overline{\pi^{p\ddagger}} \sqtimes 
\pi^{r}$, and that the operators
${\cal{P}}^{pr}_{ij}$ are the basis vectors of the carrier space 
${\cal{L}}^{pr}$ of this coaction.

Taken with (\ref{eq:V37h}), this shows that the twisted irreducible tensor 
operators
${\widetilde{Q}}^{q}_{j}$ of the definition (\ref{eq:V37e}) exist only if 
$\pi^{q}$ is contained in
the reduction of 
$\overline{\pi^{p\ddagger}} \sqtimes \pi^{r}$. As 
$\overline{\pi^{p\ddagger}}$ and $\bar{\pi^{p}}$
 are equivalent$^{7-12}$, this implies that ${\widetilde{Q}}^{q}_{j}$ exists 
only if
$n^{q}_{\bar{p}r} > 0$. But $n^{q}_{\bar{p} r} = n^{r}_{pq}$ (c.f. (I.125)), so
${\widetilde{Q}}^{q}_{j}$ exists only if $n^{r}_{pq} > 0$. (These 
observations are confirmed by the
explicit expressions for the twisted irreducible tensor operators given in 
(\ref{eq:V45h}) below and
by the Wigner-Eckart theorem of (\ref{eq:V57}) below).  

\subsubsection{Explicit expressions for the irreducible tensor operators}

If $n^{r}_{qp} > 0$ there exist $n^{r}_{qp}$ linearly independent 
{\em{ordinary}} irreducible tensor
operators that satisfy (\ref{eq:V20}). These are given by
\begin{equation} Q^{q,\alpha}_{j} = \sum_{i=1}^{d_{p}} 
\sum_{\ell=1}^{d_{r}} \left( \begin{array}{cc}
r & \bar{p}\\ \ell &  i \end{array} \right| \left. \begin{array}{ccc}q&,& 
\alpha\\j & &  
\end{array}\right) {\cal{P}}^{pr}_{i\ell} ,
\label{eq:V45d} \end{equation}   
for $\alpha = 1, 2, \ldots, n^{r}_{qp}$ and $j = 1, 2, \ldots, d_{q}$. Here 
the label $\bar{p}$ in
the Clebsch-Gordan coefficients relates to the corepresentation 
$\bar{\pi}^{p}$ that is conjugate to
$\pi^{p}$. As noted in the previous subsubsection, $n^{q}_{r\bar{p}} = 
n^{r}_{qp}$ (c.f. (I.125)). 

The proof of (\ref{eq:V45d}) is as follows. The analogue
of (I.130) for the tensor product $\pi^{r} \sqtimes \bar{\pi}^{p}$ is
\begin{equation} \sum_{\ell=1}^{d_{r}} \sum_{i=1}^{d_{p}} (\pi^{r} 
\sqtimes \bar{\pi}^{p})_{nm,\ell
i}\left(
\begin{array}{cc} r & \bar{p}\\ \ell &  i \end{array} \right| \left. 
\begin{array}{ccc}q&,&
\alpha\\j & &  
\end{array}\right) =  \sum_{k=1}^{d_{q}} \left(
\begin{array}{cc} r & \bar{p}\\ n &  m \end{array} \right| \left. 
\begin{array}{ccc}q&,& \alpha\\k &
&   \end{array}\right) \pi^{q}_{kj}
\label{eq:V45f} \end{equation} 
for $m = 1, 2, \ldots, d_{p}$, $j = 1, 2, \ldots, d_{q}$, $n = 1,
2, \ldots, d_{r}$, and $\alpha = 1, 2, \ldots, n_{qp}^{r}$. However, by 
(\ref{eq:V30}),
(\ref{eq:V45a}), and (\ref{eq:V45d}),
\[ \pi_{{\cal{L}}^{pr}}(Q^{q,\alpha}_{j}) = \sum_{i,m=1}^{d_{p}} 
\sum_{\ell
,n=1}^{d_{r}} \left(
\begin{array}{cc} r & \bar{p}\\ \ell &  i \end{array} \right| \left. 
\begin{array}{ccc}q&,&
\alpha\\j & &   \end{array}\right) \;{\cal{P}}^{pr}_{mn}\otimes (\pi^{r} 
\sqtimes
\bar{\pi}^{p})_{nm,\ell i}\;. \]
On applying (\ref{eq:V45f}) and (\ref{eq:V45d}), this reduces to 
\[ \pi_{{\cal{L}}^{pr}}(Q^{q,\alpha}_{j}) = \sum_{k=1}^{d_{q}} 
Q^{q,\alpha}_{k} \otimes
\pi^{q}_{kj}\] (for all $j = 1, 2, \ldots, d_{q}$ and $\alpha = 1, 2, \ldots,
n^{r}_{qp}$). That is, the operators $Q^{q,\alpha}_{j}$ defined in 
(\ref{eq:V45d}) satisfy
(\ref{eq:V34}), which is equivalent to (\ref{eq:V20}). 
 
Similarly, if $n^{r}_{pq} > 0$ there exist $n^{r}_{pq}$ linearly independent 
{\em{twisted}}
irreducible tensor operators that satisfy (\ref{eq:V20}). These are given by
\begin{equation} {\widetilde{Q}}^{q,\alpha}_{j} = \sum_{i=1}^{d_{p}} 
\sum_{\ell=1}^{d_{r}} \left(
\begin{array}{cc} \overline{p^{\ddagger}} & r\\ i &  \ell \end{array} 
\right| \left.
\begin{array}{ccc}q&,&
\alpha\\j & &  
\end{array}\right) {\cal{P}}^{pr}_{i\ell} ,
\label{eq:V45h} \end{equation}   
for $\alpha = 1, 2, \ldots, n^{r}_{pq}$ and $j = 1, 2, \ldots, d_{q}$. Here 
the label
$\overline{p^{\ddagger}}$ in the Clebsch-Gordan coefficients relates to 
the corepresentation
$\overline{\pi^{p\ddagger}}$ that  is conjugate to the corepresentation 
$\pi^{p\ddagger}$ that is
itself doubly contragredient to $\pi^{p}$. It should be  noted that 
$n^{q}_{\bar{p} r} =
n^{r}_{pq} = n^{q}_{\overline{p^{\ddagger}}r}$ (c.f. (I.125)). 

The proof of (\ref{eq:V45h}) is similar to that of (\ref{eq:V45d}), but uses
the right coaction $\widetilde{\pi}_{{\cal{L}}^{pr}}$ in the place of 
$\pi_{{\cal{L}}^{pr}}$.
Also needed is the relation
\[   (\overline{\pi^{p\ddagger}})_{ij} = \sum_{k,\ell=1}^{d_{p}} 
\overline{F^{p}_{ik}}\;
(\bar{\pi}^{p})_{k\ell} \;
\overline{(({\bf{F}}^{p})^{-1})_{\ell j}}\;,\]
which follows from (\ref{eq:R148}), and the corresponding relation
\[ \left(
\begin{array}{cc} \bar{p} & r\\ i &  \ell \end{array} \right| \left.
\begin{array}{ccc}q&,&
\alpha\\j & &  
\end{array}\right) = \sum_{k=1}^{d_{p}}\overline{(({\bf{F}}^{p})^{-
1})_{ik}}\;\left(
\begin{array}{cc} \overline{p^{\ddagger}} & r\\ k &  \ell \end{array} 
\right| \left.
\begin{array}{ccc}q&,&
\alpha\\j & &  
\end{array}\right)\;.\]

\section{Theorems of the Wigner-Eckart type}

If $\pi^{p}$, $\pi^{q}$, and $\pi^{r}$ are unitary irreducible 
corepresentations of
$\cal{A}$ of dimensions $d_{p}$, $d_{q}$, and $d_{r}$ respectively, 
$v^{p}_{j}$ and
$v^{r}_{\ell}$ are basis vectors belonging to the carrier spaces $V^{p}$ and 
$V^{r}$ of $\pi^{p}$ and
$\pi^{r}$ respectively, and
$Q^{q}_{k}$ is an {\em{ordinary}} irreducible tensor operator belonging to 
$\pi^{q}$ (as defined in
(\ref{eq:V20})), then
\begin{eqnarray} \langle v^{r}_{\ell},Q^{q}_{k}(v^{p}_{j}) \rangle =
\sum_{\alpha=1}^{n_{qp}^{r}} \left(\begin{array}{ccc} r&,&
\alpha\\\ell & &    
\end{array} \right| \left. \begin{array}{cc}q & p\\ k & j   
\end{array}\right)  (r \mid Q^{q} \mid
p)_{\alpha}\; ,
\label{eq:V53} \end{eqnarray}
for all $j = 1,2, \ldots, d_{p}$, all $k = 1,2,
\ldots, d_{q}$, and all $\ell = 1,2, \ldots, d_{r}$. Here the {\em{reduced 
matrix elements}} $(r \mid
Q^{q} \mid p)_{\alpha}$ are given by
\begin{eqnarray} \begin{array}{l}(r \mid Q^{q} \mid p)_{\alpha}  \\  
=\sum_{s=1}^{d_{p}}
\sum_{t=1}^{d_{q}}
\sum_{u,v=1}^{d_{r}} \langle v^{r}_{u},Q^{q}_{t}(v^{p}_{s}) \rangle 
\left( \begin{array}{cc} q & p\\ t &  s \end{array} \right| \left. 
\begin{array}{ccc}r&,& \alpha\\v
 & &  \end{array}\right) \;
\{(({\bf{F}}^{r})^{-1})_{vu}/tr(({\bf{F}}^{r})^{-1})\}\end{array}
\label{eq:V54} \end{eqnarray} for $\alpha = 1, 2, \ldots, {n_{qp}^{r}}$, 
and $\langle \; ,\;
\rangle$ denotes the inner product of $V^{r}$. Here ${\bf{F}}^{r}$ is the
matrix defined in (\ref{eq:R148}).

On the other hand, if ${\widetilde{Q}}^{q}_{k}$ is a {\em{twisted}} 
irreducible tensor operator
belonging to $\pi^{q}$ (as defined in
(\ref{eq:V37e})), then
\begin{eqnarray} \langle v^{r}_{\ell},{\widetilde{Q}}^{q}_{k}(v^{p}_{j}) 
\rangle =
\sum_{\alpha=1}^{n_{pq}^{r}} \left(\begin{array}{ccc} r&,&
\alpha\\\ell & &    
\end{array} \right| \left. \begin{array}{cc}p & q\\ j & k   
\end{array}\right)  (r \mid
{\widetilde{Q}}^{q} \mid p)_{\alpha} ,
\label{eq:V57} \end{eqnarray}
for all $j = 1,2, \ldots, d_{p}$, all $k = 1,2,
\ldots, d_{q}$, and all $\ell = 1,2, \ldots, d_{r}$, where the reduced matrix 
elements $(r \mid
{\widetilde{Q}}^{q} \mid p)_{\alpha}$ are given by
\begin{eqnarray} \begin{array}{l}(r \mid {\widetilde{Q}}^{q} \mid 
p)_{\alpha}  \\ 
 =\sum_{s=1}^{d_{p}}
\sum_{t=1}^{d_{q}}
\sum_{u,v=1}^{d_{r}} \langle v^{r}_{u},{\widetilde{Q}}^{q}_{t}(v^{p}_{s}) 
\rangle 
\left( \begin{array}{cc} p & q\\ s &  t \end{array} \right| \left. 
\begin{array}{ccc}r&,& \alpha\\v
 & &  \end{array}\right)    \;
\{(({\bf{F}}^{r})^{-1})_{vu}/tr(({\bf{F}}^{r})^{-1})\}\end{array}
\label{eq:V58} \end{eqnarray} for $\alpha = 1, 2, \ldots, {n_{pq}^{r}}$.

The results (\ref{eq:V53}) and (\ref{eq:V57}) again exhibit the classic 
Wigner-Eckart
theorem behaviour, in that they show that the
$j$, $k$, and $\ell$ dependences of the inner products
$\langle v^{r}_{\ell},Q^{q}_{k}(v^{p}_{j}) \rangle$ and
$\langle v^{r}_{\ell},{\widetilde{Q}}^{q}_{k}(v^{p}_{j}) \rangle$ are 
determined only by
Clebsch-Gordan coefficients, but it should be noted that in the general 
case in which $\cal{A}$ is
non-commutative, the inner products for the {\em{ordinary}} and 
{\em{twisted}} irreducible tensor
operators involve {\em{different}} sets of Clebsch-Gordan coefficients.

The proof of (\ref{eq:V53}) is as follows. The condition for a 
corepresentation $\pi_{V}$ of
$\cal{A}$ to be unitary is that
\begin{equation} \sum_{[v]} \langle w,v_{[1]}, \rangle S(v_{[2]}) = 
\sum_{[w]} \langle w_{[1]},v
\rangle (w_{[2]})^{*} \label{eq:V47} \end{equation}
for all $v,w \in V$, the carrier space of $\pi_{V}$, where
\begin{equation} \pi_{V}(v) = \sum_{[v]} v_{[1]} \otimes v_{[2]}\;,
 \label{eq:V47a} \end{equation}
with $v_{[1]} \in V$ and $v_{[2]} \in \cal{A}$ (c.f. (I.51)). Thus with $v = 
Q^{q}_{j}(v^{p}_{i})$
and
$w = v^{r}_{\ell}$, (\ref{eq:V47}), (\ref{eq:V47a}), (\ref{eq:V35}), and 
(\ref{eq:R5}) imply that
\begin{equation} \sum_{s=1}^{d_{p}} \sum_{t=1}^{d_{q}} \langle 
v^{r}_{\ell},Q^{q}_{t}(v^{p}_{s})
\rangle \;S(M(\pi^{q}_{tj}\otimes\pi^{p}_{si})) =
\sum_{u=1}^{d_{r}}  \langle v^{r}_{u},Q^{q}_{j}(v^{p}_{i})
\rangle (\pi^{r}_{u\ell})^{*}\;.
 \label{eq:V47b} \end{equation}
But $(\pi^{r}_{u\ell})^{*} = S(\pi^{r}_{\ell u})$, so acting on both sides 
with $S^{-1}$ (which is
well defined for a compact quantum group algebra$^{7-12}$), multiplying 
through from the left by
$(\pi^{r}_{\ell k})^{*}$, summing over $\ell$, and using the relation
$\sum_{\ell} M((\pi^{r}_{\ell k})^{*} \otimes (\pi^{r}_{\ell u})) = 
\delta_{uk}1_{\cal{A}}$,
(\ref{eq:V47b}) reduces to 
\begin{equation} \sum_{\ell=1}^{d_{r}} \sum_{s=1}^{d_{p}} 
\sum_{t=1}^{d_{q}} \langle
v^{r}_{\ell},Q^{q}_{t}(v^{p}_{s})
\rangle \;((\pi^{r}_{\ell k})^{*}\pi^{q}_{tj}\otimes\pi^{p}_{si}) =
  \langle v^{r}_{k},Q^{q}_{j}(v^{p}_{i})\rangle 1_{\cal{A}}\;.
 \label{eq:V50} \end{equation}
On acting with the Haar functional $h$, and applying (\ref{eq:R164}) and 
the relation $h(1_{\cal{A}})
= 1$, (\ref{eq:V53}) follows immediately. The proof of (\ref{eq:V57}) is 
similar.

\section{Example: Irreducible tensor operators for the standard 
deformation of the
function algebra of the compact Lie group $SU(2)$}

It is particularly interesting to study the foregoing theory for the case in 
which $\cal{A}$ is the
standard deformation  of the function algebra of the  compact Lie group 
$SU(2)$, because both this
Hopf algebra $\cal{A}$ and its dual $\cal{A}^{\prime}$ have been very 
extensively investigated, the
former in the language of `compact matrix pseudogroups' and the latter as 
the deformation
$U_{q}(sl(2))$ of the universal enveloping algebra $U(sl(2))$ of the simple 
Lie algebra $sl(2)$. 

\subsection{Structure of the standard deformation of the function algebra 
of the compact Lie
group
$SU(2)$}

It is well known (c.f. Ref. 16) that the irreducible
representations of the deformation
$U_{q}(sl(2))$ (for generic q) can be labelled by a single index $j$, which 
takes values $0,
\frac{1}{2}, 1,
\frac{3}{2}, \ldots $, the irreducible representation corresponding to $j$ 
being
$(2j+1)$-dimensional, with rows and columns that may be labeled by 
indices $m^{\prime}$ and
$m$ that take values $-j, -j+1, \ldots, j-1, j$, exactly as for the simple Lie 
algebra
$sl(2)$. To each of these representations corresponds a 
{\em{co}}representation of the dual Hopf
algebra
$\cal{A}$. Consequently the labels for corepresentations of $\cal{A}$ will 
henceforth always be
denoted by $j$ (possibly with a prime or subscript attached) and the rows 
and columns of 
the corresponding matrix coefficients will be labelled by these indices $m$ 
and $m^{\prime}$
(possibly with subscripts attached). (Although $q$ was used in all the 
other sections
of this paper to indicate an irreducible representation or corepresentation, 
in this section
it will be employed to denote the standard deformation parameter).         

All of the irreducible corepresentations $\pi^{j}$ (for $j = 0, \frac{1}{2}, 1,
\frac{3}{2}, \ldots $) may be taken to be unitary, and their matrix 
coefficients
$\pi^{j}_{m^{\prime}m}$ form a basis for ${\cal{A}}$. Moreover 
{\em{every}} matrix coefficient
$\pi^{j}_{m^{\prime}m}$ for $j \geq 1$ can be written as a polynomial in 
the matrix coefficients of
$\pi^{1/2}$, while $\pi^{0}_{00} = 1_{\cal{A}}$. Let 
\begin{eqnarray} {\mbox{\boldmath${\pi}$}}^{1/2} = 
\left(\begin{array}{cc} X & U \\ V & Y
\end{array}\right)\;,\label{eq:Ex1}\end{eqnarray}
where the entries are assumed to satisfy the relations
\begin{eqnarray}\left. \begin{array}{l}XU = q^{-1}UX, \; XV = q^{-1}VX, 
\; YU = qUY,\; YV = qVY, \\
 UV = VU,\; XY - q^{-1}UV = 1_{\cal{A}},\; YX - qUV = 1_{\cal{A}}.
\end{array}\right\}\label{eq:Ex3}\end{eqnarray}
In the language of `matrix pseudogroups' the matrix coefficients 
$\pi^{j}_{m^{\prime}m}$ are called
`quantum d-functions' and are denoted by $d^{j}_{m^{\prime}m}$. The 
work of Nomura$^{15}$ then
implies that 
\begin{eqnarray}\begin{array}{l}\pi^{j}_{m^{\prime}m} =
q^{(m^{\prime}-m)(2j-m^{\prime}+m)/2}\{[j+m^{\prime}]![j-
m^{\prime}]![j+m]![j-m]!\}\\
\;\;\times \; \sum_{a} 
\frac{q^{a(2j-m^{\prime}+m-a)}X^{j+m-a}U^{m^{\prime}-m+a}V^{a}Y^{j-
m^{\prime}-a}}
{[a]![j+m-a]![m^{\prime}-m+a]![j-m^{\prime}-a]!}  
\;,\end{array}\label{eq:Ex11}\end{eqnarray}
where $[n] = (q^{n}-q^{-n})/(q-q^{-1})$ and $[n]! = [n][n-1][n-
2]\ldots[2][1]$, and where the
sum over $a$ is over all integers such that the expressions in the q-
factorials are non-negative.  
(The present quantity $q$ is actually $q^{1/2}$ in the notation of 
Nomura$^{15,16}$). Then, for
example 
\begin{eqnarray} {\mbox{\boldmath${\pi}$}}^{1} = 
\left(\begin{array}{ccc} 
X^{2} & q^{1/2}[2]^{1/2}XU & U^{2}\\
q^{1/2}[2]^{1/2}XV & XY+qUV & q^{1/2}[2]^{1/2}UY \\
V^{2} & q^{1/2}[2]^{1/2}VY & Y^{2}
\end{array}\right)\;,\label{eq:Ex9}\end{eqnarray}
and
\begin{eqnarray} {\mbox{\boldmath${\pi}$}}^{3/2} = 
\left(\begin{array}{cccc} 
X^{3} & q[3]^{1/2}X^{2}U & q[3]^{1/2}XU^{2} & U^{3}\\
q[3]^{1/2}X^{2}V & X^{2}Y+q^{2}[2]XUV & q[2]XUY+q^{2}U^{2}V & 
q[3]^{1/2}U^{2}Y \\
q[3]^{1/2}XV^{2} & q[2]XVY+q^{2}UV^{2} & XY^{2}+q^{2}[2]UVY & 
q[3]^{1/2}UY^{2} \\
V^{3} & q[3]^{1/2}V^{2}Y & q[3]^{1/2}VY^{2} & Y^{3}
\end{array}\right)\;.\label{eq:Ex10}\end{eqnarray}
The product of any two matrix coefficients can (at least in principle) be 
deduced from the
expressions.

An alternative way of getting the product of any two matrix coefficients is 
to invoke
(\ref{eq:C45qextra}), for the Clebsch-Gordan coefficients are known for 
this $\cal{A}$. Indeed
for this $\cal{A}$ the Clebsch-Gordan coefficients exhibit two simplifying 
features. Firstly, the
multiplicity is always just 1, so the index $\alpha$ in the Clebsch-Gordan 
coefficients of 
(\ref{eq:C43}) may be omitted, and secondly, the Clebsch-Gordan 
coefficients can be taken to be
purely real. As the Clebsch-Gordan series for $\pi^{j_{1}} \sqtimes 
\pi^{j_{2}}$ is the direct sum
of $\pi^{j}$ with $j = j_{1}+j_{2}, j_{1}+j_{2}-1, \ldots,|j_{1}-j_{2}|$,  
(\ref{eq:C45qextra})
reduces in this case to 
\begin{equation}  M(\pi^{j_{1}}_{m^{\prime}_{1}m_{1}} \otimes 
\pi^{j_{2}}_{m^{\prime}_{2}m_{2}}) = 
\sum_{j=|j_{1}-j_{2}|}^{j_{1}+j_{2}} \;\sum_{m^{\prime},m=-j}^{j}\; \left(
\begin{array}{cc} j_{1} & j_{2}\\ m^{\prime}_{1} &  m^{\prime}_{2} 
\end{array} \right| \left.
\begin{array}{c}j \\m^{\prime}   \end{array}\right)
\left(\begin{array}{cc} j_{1} & j_{2}\\ m_{1} &  m_{2} \end{array} \right| 
\left.
\begin{array}{c}j \\m   \end{array}\right) \pi^{j}_{m^{\prime}m} \;. 
\label{eq:Ex29} \end{equation}
Various equivalent expressions for the Clebsch-Gordan coefficients appear 
in the literature,
but the most convenient for application here is that given by 
Nomura$^{16}$, which, in the
present notation, is
\begin{eqnarray}\begin{array}{l}\left(\begin{array}{cc} j_{1} & j_{2}\\ 
m_{1} &  m_{2} \end{array} \right| \left.
\begin{array}{c}j \\m   \end{array}\right) =\\
\;\Delta(j_{1},j_{2},j) q^{\{x(j_{1})+x(j_{2})-x(j)+2(j_{1}j_{2}+j_{1}m_{2}-
j_{2}m_{1})\}/2}\\
\;\times
\{[j_{1}+m_{1}]![j_{1}-m_{1}]![j_{2}+m_{2}]![j_{2}-m_{2}]![j+m]![j-
m]![2j+1]\}^{1/2}\\
\;\times \sum_{a} \frac{(-1)^{a} q^{-a(j_{1}+j_{2}+j+1)/2}}
{[a]![j_{1}+j_{2}-j-a]![j_{1}-m_{1}-a]![j_{2}+m_{2}-a]![j-j_{2}+m_{1}+a]![j-j_{1}-
m_{2}+a]!}\;,
\end{array}\label{eq:Ex48}\end{eqnarray}
where
$\Delta(a,b,c) = \{[-a+b+c]![a-b+c]![a+b-c]!/[a+b+c+1]!\}^{1/2}$ and $x(a) = 
a(a+1)$, and
where the sum over $a$ is over all integers such that the expressions in 
the q-factorials are
non-negative. In particular
\begin{eqnarray}\begin{array}{l}\left(\begin{array}{cc} j+\frac{1}{2} & j\\ 
m+\frac{1}{2} &  -m
\end{array}
\right| \left.
\begin{array}{c}\frac{1}{2} \\\frac{1}{2}   \end{array}\right) =
(-1)^{j-m}\, q^{-j+\frac{1}{2}+\frac{3}{2}m}\,[j+m+1]^{1/2}\,
\{[2][2j]!/[2j+2]!\}^{1/2}\end{array}\;,\label{eq:Ex43}\end{eqnarray}
and
\begin{eqnarray}\begin{array}{l}\left(\begin{array}{cc} j+\frac{1}{2} & j\\ 
m-\frac{1}{2} &  -m
\end{array}
\right| \left.
\begin{array}{c}\frac{1}{2} \\-\frac{1}{2}   \end{array}\right) =
(-1)^{j-m}\, q^{\frac{1}{2}+\frac{3}{2}m}\,[j-m+1]^{1/2}\,
\{[2][2j]!/[2j+2]!\}^{1/2}\end{array}\;.\label{eq:Ex47}\end{eqnarray}

 The action of the coproduct $\Delta$ and counit $\epsilon$ of $\cal{A}$ 
on the generators of
$\cal{A}$ is given by
\[\left.\begin{array}{ll} \Delta(X) = X \otimes X + U \otimes V, & 
\Delta(Y) = V \otimes U + Y
\otimes Y,\\
\Delta(U) = X \otimes U + U \otimes Y, & \Delta(V) = V
\otimes X + Y \otimes V,
 \end{array}\right\} \]
and 
\[ \epsilon(X) = 1, \epsilon(Y) = 1, \epsilon(U) = 0, \epsilon(V) = 0.\]
Moreover, the action of the star-operation $*$ of $\cal{A}$ on the 
generators
of
$\cal{A}$ may be taken to be 
\begin{equation} X^{*} = Y, \; Y^{*} = X, \;U^{*} = -q^{-1}V, \;V^{*} = -qU
\;, \label{eq:Ex17}\end{equation}
which implies$^{15}$ that its action on any matrix coefficient is given by
\begin{equation} (\pi^{j}_{m^{\prime}m})^{*} =
(-1)^{m-m^{\prime}}q^{m-m^{\prime}}\pi^{j}_{-m^{\prime},-m}
\;. \label{eq:Ex16}\end{equation}
As $S(\pi^{j}_{m^{\prime}m}) = (\pi^{j}_{mm^{\prime}})^{*}$ (c.f. (I.52)), 
it follows that
\begin{equation} S(\pi^{j}_{m^{\prime}m}) =
(-1)^{-(m-m^{\prime})}q^{-(m-m^{\prime})}\pi^{j}_{-m,-m^{\prime}}
\;. \label{eq:Ex19}\end{equation}
Thus
\begin{equation} S^{2}(\pi^{j}_{m^{\prime}m}) =
q^{-2(m-m^{\prime})}\pi^{j}_{m^{\prime}m}
\;. \label{eq:Ex20}\end{equation} 
As the matrix coefficients of the doubly contragredient corepresentation 
$\pi^{j\ddagger}$ are
given by $(\pi^{j}_{m^{\prime}m})^{\ddagger} = 
S^{2}\pi^{j}_{m^{\prime}m}$ (c.f. (I.57)), it
follows  that the $(2j+1) \times (2j+1)$ matrix of (\ref{eq:R148}) (which 
appears in the expressions
for the reduced matrix elements (\ref{eq:V54}) and (\ref{eq:V58}) of the 
Wigner-Eckart type
theorems) is {\em{diagonal}} and that its elements are given by
\begin{equation} F^{j}_{m^{\prime}m} = \delta_{m^{\prime}m}q^{-2(j-m)}
\;. \label{eq:Ex22}\end{equation} 

\subsection{Bosonic creation and annihilation operators as irreducible 
tensor operators}

Let $b_{1}^{\dagger}, b_{1}$ and $b_{2}^{\dagger}, b_{2}$ be two pairs of 
`deformed' bosonic
creation and annihilation operators and $N_{1}$ and $N_{2}$ the 
associated number operators whose
action on the infinite-dimensional Fock space spanned  by the occupation 
number vectors
$|n_{i}\rangle$ is given (c.f. Refs. 17,18) by
\begin{eqnarray}\left. \begin{array}{l} 
b_{i}^{\dagger}|n_{i}\rangle = [n_{i}+1]^{1/2}\,|n_{i}+1\rangle\;,\\
b_{i}|n_{i}\rangle = [n_{i}]^{1/2}\,|n_{i}-1\rangle\;,\\
N_{i}|n_{i}\rangle = 
n_{i}\,|n_{i}\rangle\;,\end{array}\right\}\label{eq:A21}\end{eqnarray}
for $i = 1, 2$, where it is assumed that the vacuum state vectors 
$|0\rangle$ are such that
$b_{i}|0\rangle = 0$ for $i = 1, 2$. It is also assumed that every member 
of the set
$\{b_{1}^{\dagger},b_{1},N_{1}\}$ commutes with every member of the set
$\{b_{2}^{\dagger},b_{2},N_{2}\}$. In the deformed generalization of the 
Jordan-Schwinger realization
of $sl(2)$ (c.f. Refs. 17,18), the basis vectors of the carrier spaces of the 
irreducible
representations of
$U_{q}(sl(2))$ are given by
\begin{equation} v^{j}_{m} = |j+m,j-m\rangle\;, 
\label{eq:A36a}\end{equation} 
and these, of course, are also the basis vectors of the carrier spaces of the 
irreducible
{\em{co}}representations of $\cal{A}$.

 Then
\begin{equation} Q^{1/2}_{1/2} = b_{1}^{\dagger}q^{-\frac{1}{2}N_{2}},
\;Q^{1/2}_{-1/2} = 
b_{2}^{\dagger}q^{\frac{1}{2}N_{1}},\label{eq:A37}\end{equation}  
and
\begin{equation} Q^{1/2}_{1/2} = qb_{2}q^{\frac{1}{2}N_{1}},
\;Q^{1/2}_{-1/2} = -b_{1}q^{-\frac{1}{2}N_{2}},\label{eq:A38}\end{equation}
are two sets of pairs of {\em{ordinary}} irreducible tensor operators that 
belong to the
2-dimensional irreducible corepresentation $\pi^{1/2}$ of $\cal{A}$.

This will now be demonstrated for the the {\em{first}} pair (\ref{eq:A37}), 
starting from the
definition (\ref{eq:V20ext}), and taking $V$ to be the direct sum of all the 
carrier spaces of all
the irreducible corepresentations of $\cal{A}$ (with just one such 
irreducible corepresentation
being included from each equivalence class). Define the right coaction 
$\pi$ of $\cal{A}$ by  
\begin{equation} \pi(v^{j}_{m}) = \sum_{m^{\prime}=-j}^{j}
v^{j}_{m^{\prime}} \otimes 
\pi^{j}_{m^{\prime}m}\;,\label{eq:Ex33}\end{equation}
for all $j = 0, \frac{1}{2}, 1, \ldots$ and $m = j, j-1, \ldots, -j$.  Then in this 
case
(\ref{eq:V20ext}) becomes
\begin{equation}  ((id \otimes M) \circ (\pi \otimes
id)
\circ (Q^{1/2}_{k}
\otimes S) \circ \pi)(v^{j}_{m}) = \sum_{\ell=-1/2}^{1/2} 
Q^{1/2}_{\ell}(v^{j}_{m}) \otimes
\pi^{1/2}_{\ell k} \;,\label{eq:Ex34}
\end{equation}
for all $j = 0, \frac{1}{2}, 1, \ldots$ and $m = j, j-1, \ldots, -j$. It will now 
be shown that this
is indeed satisfied for $k = \frac{1}{2}$. (The proof for $k= -\frac{1}{2}$ is 
similar).
By (\ref{eq:A21}), (\ref{eq:Ex19}), (\ref{eq:A37}), and (\ref{eq:Ex33}), the 
left-hand side of
(\ref{eq:Ex34}) for
$k = \frac{1}{2}$ is
\begin{eqnarray}\begin{array}{l}\sum_{m^{\prime\prime} = -j-
1}^{j}\;\sum_{m^{\prime} =
-j}^{j}\;
[j+m^{\prime}+1]^{1/2}\,q^{-\frac{1}{2}(j-m^{\prime})-(m-
m^{\prime})}\,(-1)^{-(m-m^{\prime})}\\
\times \;v^{j+\frac{1}{2}}_{m^{\prime\prime}+\frac{1}{2}} \otimes
M(\pi^{j+\frac{1}{2}}_{m^{\prime\prime}+\frac{1}{2},m^{\prime}+\frac{1
}{2}}
\otimes \pi^{j}_{-m,-
m^{\prime}})\;.\end{array}\label{eq:Ex37}\end{eqnarray}
Similarly, by (\ref{eq:A21}), (\ref{eq:A37}), and (\ref{eq:Ex33}), the right-
hand
side of (\ref{eq:Ex34}) for
$k = \frac{1}{2}$ is
\begin{equation}
[j+m+1]^{1/2}\,q^{-\frac{1}{2}(j-m)}\,\;v^{j+\frac{1}{2}}_{m+\frac{1}{2}} 
\otimes
\pi^{\frac{1}{2}}_{\frac{1}{2},\frac{1}{2}}
+[j-m+1]^{1/2}\,q^{\frac{1}{2}(j+m)}\,\;v^{j+\frac{1}{2}}_{m-\frac{1}{2}} 
\otimes
\pi^{\frac{1}{2}}_{-\frac{1}{2},\frac{1}{2}}\;,\label{eq:Ex38}\end{equation}
so it remains to show that (\ref{eq:Ex37}) reduces to (\ref{eq:Ex38}). 
However, by (\ref{eq:Ex29})
and (\ref{eq:Ex43}), (\ref{eq:Ex37}) reduces to 
\begin{eqnarray}\begin{array}{l}\sum_{m^{\prime\prime} = -j-
1}^{j}\;\sum_{m^{\prime} =
-j}^{j}\;\sum_{j^{\prime} =\frac{1}{2}}^{2j+\frac{1}{2}}\;
\{[2][2j]!/[2j+2]!\}^{1/2}\,q^{(\frac{1}{2}j-m-\frac{1}{2})}\,(-1)^{(j-m)}\\
\times\,\left(\begin{array}{cc} j+\frac{1}{2} & j\\ m^{\prime}+\frac{1}{2} 
&  -m^{\prime} \end{array}
\right|\left.\begin{array}{c}\frac{1}{2} \\\frac{1}{2} \end{array}\right)  
\left(\begin{array}{cc} j+\frac{1}{2} & j\\ m^{\prime\prime}+\frac{1}{2} 
&  -m \end{array}
\right|\left.\begin{array}{c}j^{\prime} \\m^{\prime\prime}+\frac{1}{2}-
m \end{array}\right)\\
\times\;\left(\begin{array}{cc} j+\frac{1}{2} & j\\ m^{\prime}+\frac{1}{2} 
&  -m^{\prime} \end{array}
\right|\left.\begin{array}{c}j^{\prime} \\\frac{1}{2} \end{array}\right) 
\;v^{j+\frac{1}{2}}_{m^{\prime\prime}+\frac{1}{2}}
\otimes\pi^{j^{\prime}}_{m^{\prime\prime}+\frac{1}{2}-
m,\frac{1}{2}}\;.\end{array}\label{eq:Ex43a}\end{eqnarray}
On invoking the Clebsch-Gordan orthogonality relation
\begin{eqnarray}\begin{array}{l}\sum_{m^{\prime} =-j}^{j}\;
\;\left(\begin{array}{cc} j+\frac{1}{2} & j\\ m^{\prime}+\frac{1}{2} &  
-m^{\prime} \end{array}
\right|
\left.\begin{array}{c}\frac{1}{2} \\\frac{1}{2} \end{array}\right)
\left(\begin{array}{cc} j+\frac{1}{2} & j\\ m^{\prime}+\frac{1}{2} &  
-m^{\prime} \end{array}
\right|\left.\begin{array}{c}j^{\prime} \\\frac{1}{2} \end{array}\right) 
= \left\{\begin{array}{l}1 {\mbox{, if}} \;j^{\prime} = \frac{1}{2}\\
0 {\mbox{, if}}\; j^{\prime} \neq \frac{1}{2}
\end{array}\right.\;,\end{array}\label{eq:Ex44}\end{eqnarray}
(\ref{eq:Ex43a}) (and hence (\ref{eq:Ex37})) reduces to 
\begin{eqnarray}\begin{array}{l}\sum_{m^{\prime\prime} = -j-1}^{j}\;
\{[2][2j]!/[2j+2]!\}^{1/2}\,q^{(\frac{1}{2}j-m-\frac{1}{2})}\,(-1)^{(j-m)}\\
\times \;\left(\begin{array}{cc} j+\frac{1}{2} & j\\ 
m^{\prime\prime}+\frac{1}{2} &  -m \end{array}
\right|\left.\begin{array}{c}\frac{1}{2} \\m^{\prime\prime}+\frac{1}{2}-
m \end{array}\right)
 \;v^{j+\frac{1}{2}}_{m^{\prime\prime}+\frac{1}{2}}
\otimes\pi^{j^{\prime}}_{m^{\prime\prime}+\frac{1}{2}-
m,\frac{1}{2}}\;.\end{array}\label{eq:Ex45}\end{eqnarray}
However, the remaining Clebsch-Gordan coefficients are zero if 
$m^{\prime\prime}+\frac{1}{2}-m >
\frac{1}{2}$, i.e. if $m^{\prime\prime} > m$, and are zero if  
$m^{\prime\prime}+\frac{1}{2}-m <
-\frac{1}{2}$, i.e. if $m^{\prime\prime} < m-1$, so these Clebsch-Gordan 
coefficients are non-zero
only for $m^{\prime\prime} = m, m-1$. Thus (\ref{eq:Ex45}) (and hence 
(\ref{eq:Ex37})) becomes
\begin{eqnarray*}\begin{array}{l}
\{[2][2j]!/[2j+2]!\}^{1/2}\,q^{(\frac{1}{2}j-m-\frac{1}{2})}\,(-1)^{(j-m)}\\
\times \left\{\begin{array}{c}\left(\begin{array}{cc} j+\frac{1}{2} & j\\ m-\frac{1}{2} &  -m
\end{array}
\right|\left.\begin{array}{c}\frac{1}{2} \\-\frac{1}{2} \end{array}\right)
 \;v^{j+\frac{1}{2}}_{m-\frac{1}{2}}
\otimes\pi^{\frac{1}{2}}_{-\frac{1}{2},\frac{1}{2}} +
\left(\begin{array}{cc} j+\frac{1}{2} & j\\ m+\frac{1}{2} &  -m \end{array}
\right|\left.\begin{array}{c}\frac{1}{2} \\\frac{1}{2} \end{array}\right)
 \;v^{j+\frac{1}{2}}_{m+\frac{1}{2}}
\otimes\pi^{\frac{1}{2}}_{\frac{1}{2},\frac{1}{2}}\end{array}\right\},\end{array}
\end{eqnarray*}
which, by (\ref{eq:Ex43}) and (\ref{eq:Ex47}), reduces to (\ref{eq:Ex38}).
 Similarly
\begin{equation} {\widetilde{Q}}^{1/2}_{1/2} = 
b_{1}^{\dagger}q^{\frac{1}{2}N_{2}},
\;{\widetilde{Q}}^{1/2}_{-1/2} = b_{2}^{\dagger}q^{-
\frac{1}{2}N_{1}},\label{eq:A39}\end{equation}  
and
\begin{equation} {\widetilde{Q}}^{1/2}_{1/2} = q^{-1}b_{2}q^{-
\frac{1}{2}N_{1}},
\;{\widetilde{Q}}^{1/2}_{-1/2} = 
-b_{1}q^{\frac{1}{2}N_{2}},\label{eq:A40}\end{equation}
are two sets of pairs of {\em{twisted}} irreducible tensor operators  
belonging to the
2-dimensional irreducible corepresentation $\pi^{1/2}$ of $\cal{A}$. (This 
is easily deduced from
(\ref{eq:Ex37}) and (\ref{eq:Ex38}), because in the special case of this 
algebra $\cal{A}$,
(\ref{eq:Ex3}) and (\ref{eq:Ex19}) imply that the substitutions $M 
\rightarrow M \circ \sigma$ and
$S \rightarrow S^{-1}$ merely correspond to the replacement of $q$ by 
$q^{-1}$).

It has been observed previously by Biedenharn and Tarlini$^{19}$ that 
(\ref{eq:Ex37})
provide a pair of irreducible tensor operators for the 2-dimensional 
irreducible
representation of $U_{q}(sl(2))$, their argument essentially using 
(\ref{eq:V13ext})
and the generalized Jordan-Schwinger realization of the generators of 
$U_{q}(sl(2))$, together with
various identities involving the creation and annihilation operators. The 
object of the above
analysis in this subsection is to give an explicit demonstration of the
applicability of the {\em{new}} definitions (\ref{eq:V13ext}) and 
(\ref{eq:V20exta}) for $\cal{A}$,
which, of course, apply not merely to this example but to {\em{any}} 
compact quantum group algebra.

\acknowledgements

This work was started while the author was on research leave at the 
Sektion Physik of the
Ludwig-Maximilians-Universit\"{a}t, M\"{u}nchen. The author is grateful 
to Professor J. Wess and his
colleagues for their hospitality and to the Deutscher Akademischer 
Austauschdienst for its financial
support.

\appendix

\section{Introduction}

The purpose of this Appendix is to {\em{motivate}} the definitions that are 
given in the main body of
the paper for the irreducible tensor operators. This will be done by
considering the simple special case in which the Hopf algebra $\cal{A}$ is 
the set of functions
defined on a {\em{finite}} group $\cal{G}$ of order $g$, so that the dual 
$\cal{A}^{\prime}$ of
$\cal{A}$ is the group algebra of $\cal{G}$. Of course, as $\cal{A}$ is 
commutative in this special
case, the resulting expressions are to some extent ambiguous, in that in 
this special case $M$ is
indistinguishable from $M \circ \sigma$ and $S$ is indistinguishable 
from
$S^{-1}$. The demonstration of the correctness, consistency, and 
usefulness of the definitions that
are actually employed for the {\em{general}} case are the subject matter of 
the self-contained
arguments of the main body of this paper.         

A summary of the basic facts concerning the
relationship of $\cal{A}$ and $\cal{A}^{\prime}$ may be found in the 
Introduction to the Appendix
of Paper I. 

\section{Motivation for definitions of irreducible tensor operators}

The starting point of the present argument is
(\ref{eq:Vintro2}), which of course also
applies to finite groups, and which may be
rewritten as
\begin{equation}\hat{\pi}^{\prime r}(x)\,Q^{q}_{j}\,\hat{\pi}^{\prime 
p}(x^{-1}) =
\sum_{k=1}^{d_{q}}
\Gamma^{q}(x)_{kj}Q^{q}_{k}\label{eq:V11}\end{equation} 
for all $x \in \cal{G}$ and all $j = 1,
2, \ldots, d_{q}$. Here the operators $\hat{\pi}^{\prime p}(x)$ are defined 
by
\begin{equation}\hat{\pi}^{\prime p}(x)(v^{p}_{j}) =
\sum_{k=1}^{d_{p}}
\Gamma^{p}(x)_{kj}v^{p}_{k}\;,\label{eq:V11ext}\end{equation}  
and are related to the corresponding left action $\pi^{\prime p}$ of 
$\cal{A}^{\prime}$ (a mapping of
the carrier space
$V^{p}$ into $V^{p} \otimes \cal{A}^{\prime}$) by the prescription
\begin{equation}\hat{\pi}^{\prime p}(x)(v^{p}_{j}) =
\pi^{\prime p}(x \otimes v^{p}_{j})\;.\label{eq:V21}\end{equation}
As here $S_{\cal{A}^{\prime}}(x) = x^{-1}$ and 
$\Delta_{\cal{A}^{\prime}}(x) = x \otimes x$,
(\ref{eq:V11}) can be rewritten in purely Hopf algebraic terms (for 
$\cal{A}^{\prime}$) as
\begin{equation} (\pi^{\prime r} \circ (id \otimes Q^{q}_{j}) \circ (id 
\otimes \pi^{\prime
p}) \circ (id \otimes S_{\cal{A}^{\prime}} \otimes id) \circ 
(\Delta_{\cal{A}^{\prime}} \otimes
id))(x \otimes v^{p}) =
\sum_{k=1}^{d_{q}}
\Gamma^{q}(x)_{kj}Q^{q}_{k}(v^{p})\label{eq:V13}\end{equation} 
for all $x \in \cal{A}^{\prime}$, all $v^{p} \in V^{p}$, and all $j = 1,
2, \ldots, d_{q}$. (Here the $\Gamma^{q}(x)_{kj}$ are now matrix 
elements of the irreducible
representation  $\pi^{\prime q}$ of the Hopf algebra $\cal{A}^{\prime}$).

This condition can be recast entirely in terms of quantities defined for the 
Hopf algebra
$\cal{A}$ in the following way. As noted in equation (I.A7), the 
relationship between a right
coaction
$\pi_{V}$ of 
$\cal{A}$ and the corresponding left action $\pi_{V}^{\prime}$ of 
$\cal{A}^{\prime}$ (with the
same carrier space $V$) is
\begin{equation} \pi^{\prime}_{V}(a^{\prime} \otimes v) = (\MVC \circ 
(id \otimes ev) \circ (\sigma
\otimes id) \circ (id \otimes \pi_{V}))(a^{\prime} \otimes v) 
\label{eq:R14} \end{equation} for all $a^{\prime} \in \cal{A}^{\prime}$ 
and all $v \in V$, where the
evaluation map  $ev$ (from $\cal{A}^{\prime} \otimes
\cal{A}$ to $\C$) is defined (c.f. (I.41)) by
\begin{equation} ev(a^{\prime} \otimes a) = \langle a^{\prime},a \rangle  
\label{eq:extra10} \end{equation} for all $a^{\prime} \in 
\cal{A}^{\prime}$ and all $a \in \cal{A}$.
On applying this twice (once with with $\pi^{\prime}_{V} = \pi^{\prime 
r}$ and
once with $\pi^{\prime}_{V} = \pi^{\prime p}$), the left-hand side of 
(\ref{eq:V13}) becomes
\begin{eqnarray} \begin{array}{l} (\MVrC \circ (id \otimes ev) \circ 
(\sigma \otimes id) \circ (id
\otimes \pi^{r})
\circ (id \otimes Q^{q}_{j}) \circ (id \otimes \pi^{p}) \circ (id \otimes 
\MVpC)\\ \circ (id \otimes
id \otimes ev) \circ (id \otimes \sigma \otimes id) \circ (id \otimes 
S_{\cal{A}^{\prime}}\otimes
\pi^{p}) \circ (\Delta_{\cal{A}^{\prime}} \otimes id)) (x \otimes v^{p})\;. 
\end{array}
\label{eq:V16} \end{eqnarray}  
As (I.A2) and (I.A3) can be rewritten as 
\[ (\MC \circ (ev \otimes ev) \circ (id \otimes \sigma \otimes
id) \circ (\Delta_{\cal{A}^{\prime}} \otimes id))(a^{\prime} \otimes a 
\otimes b) = (ev \circ (id
\otimes M))(a^{\prime} \otimes a \otimes b)\]  
for all $a,b \in \cal{A}$ and all $a^{\prime} \in \cal{A}^{\prime}$, and
\[ (ev \circ (S_{\cal{A}^{\prime}} \otimes id))(a^{\prime}
\otimes a) = (ev \circ (id \otimes S))(a^{\prime} \otimes a) \]
for all $a \in \cal{A}$ and all $a^{\prime} \in \cal{A}^{\prime}$, 
(\ref{eq:V16}) can be
re-expressed as 
\begin{eqnarray} \begin{array}{l} (\MVrC \circ (id \otimes ev) \circ 
(\sigma \otimes id) \circ
(id \otimes id \otimes M) \circ (id \otimes \pi^{r} \otimes id) \circ (id 
\otimes Q^{q}_{j} \otimes
S) \\
\circ (id \otimes \pi^{p}))(x \otimes v^{p})\;.
\end{array}
\label{eq:V17} \end{eqnarray}  
However, the right-hand side of (\ref{eq:V13}) can be rewritten using 
(I.A10) as
$\sum_{k=1}^{d_{q}}
\langle x , \pi^{q}_{kj} \rangle Q^{q}_{k}(v^{p})$, and hence as 
\begin{equation} \sum_{k=1}^{d_{q}}(\MVrC \circ (id
\otimes ev)
\circ (\sigma \otimes id) \circ (id \otimes Q^{q}_{k} \otimes id))(x 
\otimes v^{p}\otimes
\pi^{q}_{kj}) \;.\label{eq:V18}\end{equation} 
On equating (\ref{eq:V17}) and (\ref{eq:V18}), as the first three terms are 
common to both
expressions, they
 can be removed. The remaining terms act simply as the identity on the 
factor $x$, so on removing
this now trivial effect on $x$, it follows that (\ref{eq:V13}) is equivalent to
\[  ((id \otimes M) \circ (\pi^{r} \otimes id) \circ (Q^{q}_{j}
\otimes S) \circ \pi^{p})(v^{p}) = \sum_{k=1}^{d_{q}} Q^{q}_{k}(v^{p}) 
\otimes
\pi^{q}_{kj}
\]  for all $v^{p} \in V^{p}$ and all $j = 1, 2, \ldots, d_{q}$.
As this involves {\em{only}} quantities defined for $\cal{A}$, it provides 
the desired criterion
(\ref{eq:V20}).  

Now consider the situation in which $V$ is a vector space that is a direct 
sum of carrier
spaces of unitary irreducible corepresentations of $\cal{A}$ and which 
contains at least
$V^{p}\oplus V^{r}$. Let $\pi$ be the mapping of $V$ into $V \otimes 
\cal{A}$ that coincides
with $\pi^{p}$ on $V^{p}$ and with $\pi^{r}$ on $V^{r}$, and which acts 
similarly on any other
carrier spaces that might be contained in $V$. Of course $V$ is also is the  
direct sum of carrier
spaces of unitary irreducible representations of $\cal{A}^{\prime}$. Then 
the generalization of
(\ref{eq:V13}) to this situation is
\begin{equation} (\pi^{\prime} \circ (id \otimes Q^{q}_{j}) \circ (id 
\otimes \pi^{\prime}) \circ
(id \otimes S_{\cal{A}^{\prime}} \otimes id) \circ 
(\Delta_{\cal{A}^{\prime}} \otimes id))(x \otimes
v) =
\sum_{k=1}^{d_{q}}
\Gamma^{q}(x)_{kj}Q^{q}_{k}(v)\label{eq:V13ext}\end{equation} 
for all $x \in \cal{A}^{\prime}$, all $v \in V$, and all $j = 1, 2, \ldots, 
d_{q}$. (Here
$\pi^{\prime}$ is the mapping of $V \otimes \cal{A}^{\prime}$ into $V$ 
that
coincides with
$\pi^{\prime p}$ on
$V^{p}$ and with $\pi^{\prime r}$ on $V^{r}$, and which acts similarly 
on any other carrier
spaces that might be contained in $V$). The generalization of 
(\ref{eq:V20}) to this situation is
obviously 
\begin{equation}  ((id \otimes M) \circ (\pi \otimes id) \circ (Q^{q}_{j}
\otimes S) \circ \pi)(v) = \sum_{k=1}^{d_{q}} Q^{q}_{k}(v) \otimes
\pi^{q}_{kj}
\label{eq:V20exta} \end{equation}  for all $v \in V$ and all $j = 1, 2, 
\ldots, d_{q}$.

Because $M$ is indistinguishable from $M \circ \sigma$ and $S$ is 
indistinguishable from
$S^{-1}$ in the situation being considered here, the above arguments 
would equally well apply with
each of the following 3 substitutions:
\begin{enumerate}
\item replace $M$ by $M \circ \sigma$, but leave $S$ unchanged; 
\item leave $M$ unchanged, but replace $S$ by $S^{-1}$;
\item replace $M$ by $M \circ \sigma$ {\em{and}} replace $S$ by $S^{-
1}$.
\end{enumerate} However, in the general case in which $\cal{A}$ is non-
commutative,  the
possibilities (1) and (2) are {\em{excluded}} because with them, and in the 
situation discussed in
the previous paragraph, the identity operator would not be an irreducible 
tensor operator belonging
to the identity corepresentation. With the substitution (3), (\ref{eq:V20}) 
changes into
(\ref{eq:V37e}), which is the defining condition for a {\em{twisted}} 
irreducible tensor operator
${\widetilde{Q}}^{q}_{j}$. (Of course the corresponding substitutions for 
$\cal{A}^{\prime}$ are  
$\Delta_{\cal{A}^{\prime}} \rightarrow \sigma \circ 
\Delta_{\cal{A}^{\prime}}$ and
$S_{\cal{A}^{\prime}} \rightarrow (S_{\cal{A}}^{\prime})^{-1}$, so that 
the analogues of
(\ref{eq:V13}) and (\ref{eq:V13ext}) are
\begin{eqnarray} \begin{array}{l} (\pi^{\prime r} \circ (id \otimes 
{\widetilde{Q}}^{q}_{j}) \circ
(id \otimes
\pi^{\prime p}) \circ (id \otimes (S_{\cal{A}}^{\prime})^{-1} \otimes id) 
\circ
(\sigma \otimes id) \circ (\Delta_{\cal{A}^{\prime}}
\otimes id))(x \otimes v^{p})  \\
\;\; =\sum_{k=1}^{d_{q}}
\Gamma^{q}(x)_{kj}{\widetilde{Q}}^{q}_{k}(v^{p})\end{array}\label{eq:V1
3tw}\end{eqnarray} 
(for all $x \in \cal{A}^{\prime}$, all $v^{p} \in V^{p}$, and all $j = 1,
2, \ldots, d_{q}$), and  
\begin{eqnarray} \begin{array}{l} (\pi^{\prime} \circ (id \otimes 
{\widetilde{Q}}^{q}_{j}) \circ
(id \otimes
\pi^{\prime}) \circ (id \otimes (S_{\cal{A}}^{\prime})^{-1} \otimes id) 
\circ (\sigma \otimes id)
\circ(\Delta_{\cal{A}^{\prime}} \otimes id))(x
\otimes v) \\
\;\; =\sum_{k=1}^{d_{q}}
\Gamma^{q}(x)_{kj}{\widetilde{Q}}^{q}_{k}(v) 
\end{array}\label{eq:V13twext}\end{eqnarray} 
(for all $x \in \cal{A}^{\prime}$, all $v \in V$, and all $j = 1, 2, \ldots, 
d_{q}$).  

It is worth noting that (\ref{eq:V13}), (\ref{eq:V13ext}), (\ref{eq:V13tw}), 
and (\ref{eq:V13twext})
provide the appropriate definitions for irreducible tensor operators not 
merely for the context in
which they have been derived here (i.e. for the case in which 
$\cal{A}^{\prime}$ is the group algebra
of a finite group $\cal{G}$), but also for the case in which 
$\cal{A}^{\prime}$ is the universal
enveloping algebra $U(\cal{L})$ of a Lie algebra $\cal{L}$ (with 
$S_{\cal{A}^{\prime}}(a) = -a$
and  $\Delta_{\cal{A}^{\prime}}(a) = a \otimes 1 + 1 \otimes a$ for all $a 
\in \cal{L}$), and for
deformations of such universal enveloping algebras. (Of course in the case 
$\cal{A}^{\prime} =
U(\cal{L})$, the criteria (\ref{eq:V13}) and (\ref{eq:V13tw}) coincide and 
the
criteria (\ref{eq:V13ext}) and (\ref{eq:V13twext}) also coincide, but this 
will not be true for
deformations of  $U(\cal{L})$).

\section{Derivation of the right coactions}

Consideration of (\ref{eq:Vext1}) suggests that one first defines an 
operator
$\pi^{\prime}_{{\cal{L}}^{pr}}$ to be the mapping of
${\cal{A}}^{\prime}
\otimes {{\cal{L}}^{pr}}$ into ${{\cal{L}}^{pr}}$ that is given by 
\begin{equation} \pi^{\prime}_{{\cal{L}}^{pr}}(x \otimes Q) = 
\hat{\pi}^{\prime
r}(x)\;Q\;\hat{\pi}^{\prime p}(x^{-1})\;,\label{eq:V25}\end{equation}
where the operators
$\hat{\pi}^{\prime p}(x)$ were defined (\ref{eq:V11ext}) (and the 
$\hat{\pi}^{\prime r}(x)$ are
defined similarly), and where $Q$ is any member of ${{\cal{L}}^{pr}}$. If 
$a^{1}, a^{2}, \dots $
form a basis for ${\cal{A}}^{\prime}$,  (\ref{eq:V25}) can be re-expressed 
in
purely Hopf algebra terms (with $x = a^{k}$) as 
\begin{eqnarray} \begin{array}{l} \pi^{\prime}_{{\cal{L}}^{pr}}(a^{k} 
\otimes Q) = \\(\widehat{M}
\circ (id
\otimes \widehat{M}) \circ (\hat{\pi}^{\prime r} \otimes id \otimes  
\hat{\pi}^{\prime
p}) \circ (id \otimes \sigma) \circ (id \otimes S_{\cal{A}^{\prime}} 
\otimes id) \circ
(\Delta_{\cal{A}^{\prime}}
\otimes id)) (a^{k} \otimes Q). \end{array} \label{eq:V26}\end{eqnarray}
(Here the operator multiplication operation $\widehat{M}$ is defined by
$  \widehat{M}(Q \otimes Q^{\prime}) = Q \circ Q^{\prime}$ for all 
$Q,Q^{\prime} \in
{{\cal{L}}^{pr}}$). It is then easily demonstrated that 
$\pi^{\prime}_{{\cal{L}}^{pr}}$ is a
{\em{left action}} of  ${\cal{A}}^{\prime}$ with carrier space 
${{\cal{L}}^{pr}}$. 

 After some algebra, it can be shown that (\ref{eq:V26}) can be rewritten 
in terms of
components as
\begin{equation} \pi^{\prime}_{{\cal{L}}^{pr}}(a^{k} \otimes Q) = 
\sum_{i,j = 1}^{d_{p}} \sum_{m,n
= 1}^{d_{r}}
\langle a^{k}\,,\, M(\pi^{r}_{mn} \otimes S(\pi^{p}_{ij})) \rangle\;
q_{ni} {\cal{P}}^{pr}_{jm}\;, \label{eq:V27}\end{equation}
where the operators ${\cal{P}}^{pr}_{bi}$ are defined in (\ref{eq:V5}) and 
the matrix elements
$q_{ja}$ are defined in (\ref{eq:V7}). The corresponding right coaction 
$\pi_{{\cal{L}}^{pr}}$ of
$\cal{A}$ is then given (c.f. (I.A6)) by
\begin{equation} \pi_{{\cal{L}}^{pr}}(Q) = \sum_{k}
 \pi^{\prime}_{{\cal{L}}^{pr}}(a^{k} \otimes Q) \otimes a_{k} 
\label{eq:V28}\end{equation}
for all $Q \in {{\cal{L}}^{pr}}$, where $a_{1}, a_{2}, \ldots $ is the dual 
basis of $\cal{A}$. Thus,
by (\ref{eq:V27}) and (\ref{eq:V28}), 
\[ \pi_{{\cal{L}}^{pr}}(Q) =
\sum_{i,j=1}^{d_{p}}\;\sum_{m,n=1}^{d_{r}}\; (q_{ni}\,{\cal{P}}^{pr}_{jm})
\otimes M(\pi^{r}_{mn} \otimes S(\pi^{p}_{ij}))\; ,\]
which is (\ref{eq:V30}).

On replacing $M$ by $M \circ \sigma$ and $S$ by $S^{-1}$, the definition 
(\ref{eq:V30}) changes
into the definition (\ref{eq:V37f}) for $\widetilde{\pi}_{{\cal{L}}^{pr}}$.

\end{document}